\documentclass[11pt, twoside, a4paper]{article}
\usepackage[english]{babel}

\usepackage{amsmath,amssymb,amsthm,epsfig}
\begin{document}

\newtheorem{exemp}{Exemplo}

\newtheorem{theorem}{Theorem} \newtheorem{corollary}{Corollary}
\newtheorem{lemma}{Lemma} \newtheorem{proposition}{Proposition}
\newtheorem{definition}{Definition}

\newcommand{\z}{\mathbb{Z}} \newcommand{\re}{\mathbb{R}}
\newcommand{\tn}{\mathbb{T}^N}
\newcommand{\rn}{\mathbb{R}^N}
\newcommand{\espaco}{[0,1]^{\mathbb{N}}}
\newcommand{\supp}{\mbox{supp}}
\newcommand{\nat}{\mathbb{N}}
\newcommand{\Rm}{{\noindent \sc Note. \ }}

\def\cqd {\,\hfill  \begin{footnotesize}$\square$\end{footnotesize}}
\def\cqdt {\hspace{5.6in} \begin{footnotesize}$\blacktriangleleft$\end{footnotesize}}

\def\appendixname{\empty}

\def\refname{References}
\def\bibname{References}
\def\chaptername{\empty}
\def\figurename{Fig.}
\def\abstractname{Abstract}

\title{  Negative Entropy, Zero temperature and stationary Markov chains on the
interval. }
\author{
A. O. Lopes*, J. Mohr* , R. R. Souza* \footnote{Instituto de
Matem\'atica, UFRGS, 91509-900 Porto Alegre, Brasil. Partially
supported by CNPq, PRONEX -- Sistemas Din\^amicos, Instituto do
Mil\^enio, and beneficiary of CAPES financial support. J. Mohr was
partially supported by CNPQ PhD and also Pos-doc scholarship.} and
Ph. Thieullen\footnote{Institut de Math\'ematiques, Universit\'e
Bordeaux 1, F-33405 Talence, France. }}
\date{\today}
\maketitle

\begin{abstract}

We consider ergodic optimization for the shift map on the modified Bernoulli
 space $\sigma:\espaco\to [0,1]^\mathbb{N}$, where $[0,1]$ is the unit closed interval, and the
 potential $A:[0,1]^\mathbb{N}\rightarrow \mathbb{R}$ considered depends on the two first
 coordinates of $[0,1]^\mathbb{N}$.
 We are interested in
finding stationary Markov probabilities $\mu_\infty$ on $
[0,1]^\mathbb{N}$ that maximize the value $ \int A d \mu,$ among all
stationary (i.e. $\sigma$-invariant) probabilities $\mu$ on
$[0,1]^\mathbb{N}$. This problem correspond in Statistical Mechanics
to the zero temperature case for the interaction described by the
potential $A$. The main purpose of this paper is to show, under the
hypothesis of uniqueness of the maximizing probability, a Large
Deviation Principle for a family of absolutely continuous Markov
probabilities $\mu_\beta$ which weakly converges to $\mu_\infty$.
The probabilities $\mu_\beta$ are obtained via an information we get
from a Perron operator and they satisfy a variational principle
similar to the pressure in Thermodynamic Formalism. As the potential
$A$ depends only on the first two coordinates, instead of  the
probability $\mu$ on $[0,1]^\mathbb{N}$,  we can consider its
projection $\nu$ on $[0,1]^2$. We look at the problem in both ways.
If $\mu_\infty$ is the maximizing probability on $[0,1]^\mathbb{N}$,
we also have that its projection $\nu_\infty$ is maximizing for $A$.
The hypothesis about stationarity on the maximization problem can
also be seen as a transhipment problem. Under the hypothesis of $A$
being $C^2$ and the twist condition, that is, $\frac{\partial^2 \, A
}{\partial x \partial y} (x,y) \neq 0$, for all $(x,y) \in [0,1]^2$,
we show the graph property of the maximizing probability $\nu$ on
$[0,1]^2$. Moreover, the graph is monotonous. An important result we
get is: the maximizing probability is unique generically in
Ma\~n\'e's sense. Finally, we exhibit a separating sub-action for
$A$.
\newline
AMS  28D05; 60J10; 37C40; 82B05
\end{abstract}

\newpage

\vspace {.8cm}


\bigskip
\section{Introduction }\label{secaoinicial}

We consider ergodic optimization \cite{Jen1}
 \cite{CG} \cite{CLT} \cite{Mo} for the shift map on the modified Bernoulli
 space $\sigma:\espaco\to [0,1]^\mathbb{N}$, where $[0,1]$ is the unit closed interval, and the
 potential $A:[0,1]^\mathbb{N}\rightarrow \mathbb{R}$ considered depends on the two first
 coordinates of $[0,1]^\mathbb{N}$.
 We are interested in
finding stationary Markov probabilities $\mu_\infty$ on $
[0,1]^\mathbb{N}$ that maximize the value $ \int A d \mu,$
among all stationary (i.e. $\sigma$-invariant) probabilities $\mu$ on
$[0,1]^\mathbb{N}$, and study properties of this maximizing measures.


 We denote by $\mathbf {x}=(x_1,x_2,....)$ a point in $\espaco$,
and we consider the shift map $\sigma:\espaco \rightarrow \espaco$ given by $\sigma((x_1,x_2,...))=(x_2,x_3,...)$.
The sigma-algebra we consider in $\espaco$ is the one generated by the cylinders.

By a stationary probability (or stationary measure) we mean a  probability  that is $\sigma$-invariant.
By a stationary Markov probability we mean a stationary probability that is obtained from an initial probability $\theta$ on $[0,1]$, and a Markovian transition Kernel $dP_x(y)=P(x,dy)$, where $\theta$ is invariant for the kernel  defined by $P$.
In the next section we will present precise definitions.

We consider a continuous potential  $A:\espaco \to \mathbb{R}$ which depends only on the two first coordinates of $\espaco$.  Therefore, we can define $\tilde A:[0,1]^2\to \mathbb{R}$, as $ \tilde A(x_1,x_2)=A(\mathbf {x})$, where   $\mathbf {x}$ is any point in $\espaco$ which has $x_1$ and $x_2$ as its two first coordinates. We will drop the symbol $\,\tilde{}\, $ and the context will show if we are considering a potential in $[0,1]^2$ or in $\espaco $.


 We are interested in finding stationary Markov probabilities
$\mu_\infty$ on the Borel sets of  $\espaco $ that maximize the value
$$ \int A(x_1,x_2) \,d \mu(\mathbf {x}) ,$$
among all stationary probabilities $\mu$ on $\espaco$.

The maximizing probabilities $\mu_\infty$, in general, 
are not positive
in all open sets on $\espaco$.

 We present an entropy penalized method
(see \cite{GV} for the case of Mather measures) designed to
approximate a maximizing probability $\mu_\infty$ by (absolutely
continuous) stationary Markov probabilities $\mu_\beta$, $\beta>0$,
obtained from ${\theta}_{\beta}(x)$ and ${P}_{\beta}(x,y)$ which are
continuous functions. The functions ${\theta}_{\beta}$ and ${P}_{\beta}$ are
obtained from the eigenfunctions and the eigenvalue   of a pair of  Perron operator (we consider
the operators $\varphi \to            {\cal L}_{\beta} \varphi (\cdot)
=\int e^{\beta\, A(x,\cdot)} \, \varphi(x)dx$ and $\varphi \to            \bar{\cal L}_{\beta} \varphi (\cdot)
=\int e^{\beta\, A(\cdot,y)} \, \varphi(y)dy$  and we use Krein-Ruthman
Theorem) in an analogous way as  the case described by F. Spitzer in
\cite{Sp} for the Bernoulli space $\Omega=\{1,2,..,d\} ^\mathbb{N}$
(see also \cite{PP}).

We will show a large deviation principle for the sequence $\{\mu_\beta\}$  which converges to $\mu_{\infty}$ when ${\beta\to \infty}$.
The large deviation principle give us important information on the rate of such convergence  \cite{DZ}.

\vspace{.2cm}

When the state space is the closed unit interval $[0,1]$, therefore,
not countable, strange properties can occur: the natural variational
problem of pressure deals with a negative entropy, namely, we have
to consider the entropy penalized concept. Negative entropies appear in a natural way when we  deal with
a continuous state space (see \cite{Ju} for mathematical results and also applications to Information Theory).
In physical problems they occur  when the spins are in a continuum space (see
for instance \cite{Lu} \cite{Cv} \cite{Ni}  \cite{RRS}
\cite{W} \cite{BBNg}).

Our result is similar to \cite{BLT} which considers the states space $S=\{1,2,..,d\}$ and \cite{GLM}
which consider the entropy penalized method for Mather measures \cite{CI} \cite{Fathi}.

In a certain extent, the problem we consider here can be analyzed just by considering probabilities $\nu$ on
$[0,1]\times[0,1]$ defined by
$$\nu( \,[a_1,a_2]\times [b_1,b_2]\,)= \int_{a_1}^{a_2} \int_{b_1}^{b_2}
 dP_{x_1}(x_2)d\theta (x_1),$$ instead of probabilities $\mu$ on $\espaco$ defined by corresponding $\theta $ and the
markovian kernel $P_x(y)$. We say that $\nu$ is the projection  of $\mu$ on $[0,1]\times[0,1]$.

From the point of view of Statistical Mechanics we are analyzing a
system of neighborhood interactions described by $A(x,y)$ at
temperature zero, where the spin $x$ takes values on $[0,1]$. This
is another point of view for the meaning of the concept of
maximizing probability for $A$. A well known example is when
$A(x,y)=x\, y$, and $x,y\in [-1,1]$  (see \cite{Th}  for
references), which can be analyzed using the methods described here
via change of coordinates. In the so called $XY$ spin model, we have
$A(x,y)=\cos(x-y)$, where $x,y \in (0, 2 \pi ]$ (see \cite{V} \cite{Pe} and  \cite{Ta} for explicit solutions). When
there is magnetic term one could consider, for instance,
$A(x,y)=\cos(x-y)+ \, l\, \cos(x)$, where $l$ is constant
\cite{RRS} \cite{A}. We show, among other things, that for
such model, given  a generic $f$ (in the sense of Ma\~n\'e
\cite{Man}), the maximizing probability for $A$ is unique. Our result seems to be related to section III b) in \cite{CG}.

Finally, another point of view for our main result: consider the cost $A:[0,1]\times[0,1]\to \mathbb{R}$, and
the problem of maximizing $\int A(x,y) \,d \nu(x,y)$, among
probabilities $\nu$ over $[0,1]\times[0,1]$ (which can be
disintegrated as $d\nu(x ,y)=d\theta (x) dP_{x}(y)$\;) with the
property of having the same marginals in the $x$ and $y$
coordinates. We refer the reader to \cite{Ra} for a broad
description of the Monge-Kantorovich mass transport problem and  the
Kantorovich-Rubinstein mass transhipment problem. We consider here a special case of such problem.
In this way we obtain a robust method (the
LDP is true) to approximate the probability $\nu_\infty$, which is
solution of the mass transhipment problem, via the entropy penalized
method.

 Under the twist
hypothesis, that is $\frac{\partial^2 \, A}{\partial x \partial y}
(x,y) \neq 0$, for all $(x,y) \in [0,1]^2,$ we show that
the probability $\nu_\infty$  on $[0,1]^2$ is supported in a graph.

The twist condition  is essential in  Aubry Theory for twist maps
\cite{Ban} \cite{Go}. It corresponds, in the  Mather Theory, to the hypothesis
of convexity of the Lagrangian \cite{Mat} \cite{CI} \cite{Fathi} \cite{Man}. It is also considered in discrete time for  optimization problems as  in \cite{Ba} \cite{Mi}. Here, several results can be
obtained without it. But, for getting results like the graph
property, it is  necessary.

In section \ref{mainres} we present some basic definitions and the main results of the paper.
In section \ref{secao_statmark} we present the induced Markov measures on $[0,1]^2$ and its relation with stationary measures on $\espaco$. In section \ref{section_maximprob} we
introduce the Perron operator, the entropy penalized concept and we
consider the associated variational problem. In section \ref{uniqueness}, under the
hypothesis of $A$ being $C^2$ and the twist condition, we show the
graph property of the maximizing probability. We also show that for the potential $A$, in
the generic sense of Ma\~n\'e (see \cite{Man} \cite{BC} \cite{CI} \cite{CLT}), the
maximizing probability on $[0,1]^2$ is unique. We get the same results for calibrated sub-actions. In section \ref{sec_ldp}, we
 present the deviation function $I$ and show the L.D.P.. In section \ref{sect_separating}, we show the monotonicity of the graph and we exhibit a separating sub-action.

All results presented here can be easily extended to Markov Chains with state space $[0,1]^2$, or, to more general potentials depending on a finite number of coordinates in $\espaco$, that is, to $A$ of the form  $A(x_1,x_2,...,x_n),$
$ A:[0,1]^n \to \mathbb{R}$.

We would like to thanks Alexandre Baraviera and Ana Ribeiro-Teixeira for references and interesting conversations on the subject
of the paper.
\bigskip

\subsection{Main results}\label{mainres}
Next we will give some definitions in order to state the main results of this work.

$\espaco$ can be endowed with the product topology, and then $\espaco$ becomes a compact metrizable topological space.
We will define a distance in $\espaco$ by
$$d(\mathbf{x},\mathbf{y})=\sum_{j\geq 1}\frac{|x_j-y_j|}{2^j}.
$$

\begin{definition}
(a) the {\em shift map} in $\espaco$ is defined as $\sigma((x_1,x_2,...)) = (x_2,x_3,...)$.
\bigskip

(b) Let $A_1, A_2, ... , A_k$ be non degenerated intervals of  $[0,1]$. We call a {\em cylinder of size  $k$ }the subset of
$\re^k$  given by   $A_1 \times A_2 \times ... \times A_{k}$, and we denote it by   $A_1....A_k$.

\bigskip

(c) Let  $\mathcal {M}_{\espaco}$ be the   set of probabilities in the
Borel sets of $[0,1]^{\mathbb{N}}$. We define the set of {\em holonomic measures} in     $\mathcal M_{\espaco} $  as
$$\mathcal M_0:=\left\{\mu\in\mathcal M_{\espaco} : \int (f(x_1)-f(x_2))\;d\mu(\mathbf{x})=0,\;\;\;\forall f\in C([0,1])\right\} \;.$$


\end{definition}
{\bf Remark}:
(i) A cylinder can also be viewed as a subset of $\espaco$: in this case, we have
$$ A_1....A_k = \left\{ \mathbf{x}\in \espaco \;: \; x_i \in A_i ,\; \forall \, 1\leq i \leq i \right\}\,.$$

(ii) For the  set of holonomic probabilities $\mathcal M_0$, we keep the terminology  used in  \cite{Gom} and  \cite{GL}. This set has been also considered in \cite{Man} and \cite{FS}.

(iii)  $\mathcal M_0$ contain all $\sigma$-invariant measures. This is a consequence of the fact that invariant measures for a transformation defined in a compact metric space can be characterized by the measures $\mu$ such that $\int f d\mu = \int (f \circ \sigma) \,d\mu$ for all continuous functions defined in $\espaco$ and taking values in $\re$. Note that 
the set of $\sigma$-invariant measures is a proper subset of $\mathcal M_0$.

\begin{definition}  A function
 $P:[0,1]\times \mathcal A\to [0,1]$ is called a {\em transition probability function} on $[0,1]$, where $\mathcal A$ is the Borel $\sigma$-algebra on $[0,1]$, if

 \vspace{.3cm}

 (i) for all $x\in[0,1]$, $P(x,\cdot)$ is a probability measure on $([0,1],\mathcal A)$,

 \vspace{.3cm}

 (ii) for all $B\in\mathcal A$, $P(\cdot, B)$ is a $\mathcal A$-measurable function from $([0,1],\mathcal A)\to [0,1]$ .
 \end{definition}

Sometimes we will use the notation $P_x(B)$ for $P(x,B)$.

Any probability $\nu$ on $ [0,1]^2$ can be disintegrated as $d\nu(x,y)= d\theta(x) dP_x(y)$, and we will denote it by $\nu=\theta P$, where
$\theta$ is a probability on $([0,1],\mathcal A)$    \cite{Dellach}, Pg 78, (70-III).

\begin{definition} \label{est}A   probability measure $\theta$ on $([0,1],\mathcal A)$ is called
{\em stationary} for a transition  $P(\cdot,\cdot)$, if
$$\theta (B)=\int P(x,B)d\theta(x)\;\;\;\;\;\;\;\mbox{ for all } B\in\mathcal A.$$

\end{definition}

Given the initial probability $\theta$  and the transition $P$, as
above, one can define a Markov process $\{X_n\}_{n\in \mathbb{N}}$
with state space $S=[0,1]$ (see \cite{AL} section 14.2  for general references on the topic). If $\theta$ is stationary for $P$, then,
one can prove that $X_n$ is a stationary stochastic process. The associated probability $\mu$ over
$[0,1]^\mathbb{N}$ is called the Markov stationary probability
defined by $\theta $ and $P$.

\begin{definition} \label{medmarkov}A probability measure   $\mu\in\mathcal
  M_{\espaco}$ will be called  a
 {\em stationary Markov measure} if there exist  $\theta$ and $P$ as in the definition
  \ref{est}, such that $\mu$ is given by

\begin{equation} \mu(A_1...A_n):=\int_{A_1...A_n}\; dP_{x_{n-1}}( x_n)...dP_{x_1}( x_2) \;d\theta(x_1) \, , \end{equation}
where  $A_1...A_n$ is a cylinder of size $n$.
\end{definition}

\bigskip

We  consider the following problem: to find measures that  maximize, over $ \mathcal M_0$, the value

$$ \int A(x_1,x_2) \,d \mu(\mathbf {x}) ,$$ which is more general than the problem of maximizing $\int A d\mu$ over  the stationary probabilities.
\bigskip

We define  
$$m=\max_{\mu\in\mathcal M_0} \left\{ \int A d\mu \right\}\,.$$
We will see that this two problems are equivalents, as   we will construct a stationary Markov measure $\mu$ such that  $m=\int A\;d\mu$.
This measure will be called a maximizing stationary Markov measure.

\begin{definition}\label{sub}
(a) A continuous function $u: [0,1] \to \mathbb{R} $ is called a
{\em calibrated forward-subaction}  if, for any $y$ we have

\begin{equation}\label{c} u(y)=\max_x [A(x,y)+ u(x)-m].\end{equation}

(b) A continuous function $u: [0,1] \to \mathbb{R} $ is called a
{\em calibrated backward-subaction} if, for any $x$ we have
\begin{equation}\label{d} u(x)=\max_y [A(x,y)+ u(y)-m].\end{equation}
\end{definition}

{\bf Remark:} If $A$ depends on all coordinates in $\espaco$,  a calibrated forward-subaction (see \cite{BLT}, but note that there they call it a strict subaction, see also \cite{GL}) is a continuous function $u:\espaco \to \mathbb{R} $ satisfying $$u(\mathbf z)=\max_{\mathbf x :\sigma(\mathbf x)=\mathbf z}[A(\mathbf x)+u (\mathbf x)-m]. $$  Hence, if $A$ depends only on the two first coordinates of $\espaco$, definition \ref{sub} is a particular case of this definition.



We denote by $C^2([0,1])$ the set of twice continuously differentiable maps from $[0,1]$ to the real line.
The main results of this paper can be summarized by the following theorems (although in the text they  will be split in several other  results):

\begin{theorem} \label{principal1} If $A$ is  $C^2$ and satisfies $\frac{\partial^2 A}{\partial x \partial y} \neq 0$, then
 there exists a generic set $\mathcal{O}$  in $C^2([0,1])$  (in Baire sense) such that:
 \medskip

 (a) for each $f \in \mathcal O$,
given $\mu, \tilde \mu\in \mathcal M_0$  two maximizing measures for $A+f$ (i.e., $m=\int (A+f)\;d\mu=\int (A+f)\;d\tilde\mu$),
then
$$\nu=\tilde \nu,$$
where $\nu$ and $\tilde \nu$ are the projections of $\mu$ and $\tilde \mu$ in the first two coordinates.

\bigskip

(b) The calibrated backward-subaction (respectively, calibrated forward-subaction) for $A+f$ is unique.
\end{theorem}

\begin{theorem}\label{principal} Let  $A:\espaco\to\re$ be a continuous  potential that depends only on the first two coordinates of  $\espaco$.
Then

 (a) There exist a measure  $\mu_{\infty}\in\mathcal M_0$ such that
   $\int A d\mu_{\infty}=m$ ,  and a sequence of stationary Markov  measures  $\mu_{\beta}$, $\beta\in\re$
  such that  $$\mu_{\beta}\rightharpoonup \mu_{\infty},$$ where $\mu_{\beta}$ is defined by  $\theta_{\beta}:[0,1]\to\re, K_{\beta}:[0,1]^2\to\re$ (see equations \eqref{theta} and \eqref{K})  as
$$ \mu_{\beta}(A_1...A_n):=\int_{A_1...A_n}
K_
{\beta}(x_{n-1},x_n)...K_{\beta}(x_1,x_2)\theta_{\beta}(x_1)dx_n...dx_1 $$
for any cylinder $ A_1...A_n$. Also $\mu_{\infty}$ is a stationary Markov measure.

\medskip

 (b) If $A$  has only one   maximizing stationary Markov  measure and there exist an unique  calibrated  forward-subaction $V$ for $A$, then the following LDP is true:
for each cylinder $D=A_1....A_k$, the following limit exists

$$\lim_{\beta\to\infty}\frac{1}{\beta}\ln
\mu_{\beta}(D)= -\inf_{\mathbf{x}\in D}I(\mathbf{x})\,.$$

where
$I:[0,1]^{\nat}\to [0,+\infty]$ is a function defined by
$$I(\mathbf{x}):= \sum_{i\geq 1}V(x_{i+1})-V(x_i)-(A-m)(x_i,x_{i+1})\,.   $$

\end{theorem}

{\bf Remark to Theorem \ref{principal}(b):} we will show, in what follows, that Theorem \ref{principal1}(a) implies  that, for $f \in \mathcal O$, the maximizing stationary Markov measure for $A+f$ is unique.




\section{ Induced stationary Markov measures}\label{secao_statmark}

In this section we consider a special class of two-dimensional measures that is closely related to the stationary measures. We will prove that the two-dimensional measure of this class that maximizes the integral of the observable $A$ can be extended to a Markov stationary measure that solves the problem of maximization of the integral of $A$ among all stationary measures.


We will denote by $\mathcal {M}_{[0,1]^2}$ the set of probabilities measures
in the Borel sets of $[0,1]^2$.
$\mathcal {M}_{[0,1]^2}$ can be endowed with the weak-$\star$ topology, where a sequence $\nu_n \to \nu$, iff, $\int f d\nu_n \to \int f d\nu$, for all continuous functions $f:[0,1]^2\to\re$. We remember that Banach-Alaoglu theorem implies that $\mathcal {M}_{[0,1]^2}$  is  a compact topological space.

\begin{definition} (a) A probability measure $\nu\in\mathcal M_{[0,1]^2}$ will be called a
 {\em induced stationary Markov measure} if its disintegration $\nu=\theta P $ is such that the probability measure $\theta$  on $([0,1],\mathcal A)$ is stationary for $P$.

  In this case for  each set $(a,b)\times(c,d) \in [0,1]^2$ we have
$$
\nu((a,b)\times(c,d)) = \int_{(a,b)} \int_{(c,d)} dP_x(y) d\theta(x)
$$

(b) We will denote by $ \mathbf{M}$ the set of induced stationary Markov measures.
\end{definition}

\begin{definition} (a) A probability measure $\nu$ will be called
an {\em induced absolutely continuous stationary Markov measure}, if
$\nu$  is in  $ \mathbf{M}$    and     can be disintegrated as $\nu=\theta K$, where
$\theta$ is an absolutely continuous measure given by a continuous density $\theta(x) dx$, and also  for each $x\in [0,1]$ the measure
$K(x,.)$ is an absolutely continuous measure given by a continuous density $K(x,y) dy$.

(b) We will denote by $\mathbf{M}_{ac}$ the set of induced absolutely continuous stationary Markov measures.
\end{definition}

We can see that the above continuous densities
$K:[0,1]^2\to[0,+\infty)$ and  $\theta:[0,1]\to[0,+\infty)$ satisfy the following equations:

\begin{equation}\label{transicao}
\int K(x,y)\, dy=1, \hspace{1cm} \forall \, x \in [0,1],
\end{equation}

\begin{equation}\label{normalizada}   \int \theta(x)\, K(x,y) \,dx dy=1, \end{equation}

\begin{equation}\label{equilibrio}  \int \theta(x)\, K(x,y) \,dx = \theta(y),\hspace{1cm} \forall \, y \in [0,1].
\end{equation}
Moreover, any pair of non-negative continuous functions satisfying the three equations above define an induced absolutely continuous stationary Markov measure.

Let  $C[0,1]$ denotes the set of continuous functions defined in $[0,1]$ and taking values on $\re$,
and  $C([0,1]^2)$ denotes the set of continuous functions defined in $[0,1]^2$ and taking values on $\re$.

\begin{lemma}\label{holon} $$\mbox{ (a) }   \mathbf{M} =\left\{\nu\in \mathcal {M}_{[0,1]^2} \; : \; \int f(x)-f(y)\;d\nu(x,y)=0
\;,\; \forall f\in C[0,1]\right\}\ ,$$

(b)  $\mathbf{M}$ is a closed set in the weak-$\star$ topology.

\end{lemma}

\textit{Proof: }(a)  Suppose that $\nu=\theta P \in \mathbf{M}$ is a induced stationary Markov measure. Remembering that $\int f d\nu$ is defined by the limit of integrals of simple functions, it is enough to show that $\int f(x) d\nu(x,y)=\int f(y) d\nu(x,y)$
 for $f=\chi_B$ where  $B$ is a Borel set. We have
$$\int_{[0,1]} \int_{[0,1]} \chi_B(x) dP_x(y) d\theta(x)=\int_{[0,1]}  \chi_B(x) \int_{[0,1]} dP_x(y) d\theta(x)=$$
$$=\int_{[0,1]}  \chi_B(x) d\theta(x)
= 
\theta(B) = \int P(x,B) d\theta(x)=$$
$$=\int_{[0,1]} \int_B dP_x(y) d\theta(x)=\int_{[0,1]} \int_{[0,1]} \chi_B(y) dP_x(y) d\theta(x)\,.$$

Now we will suppose that $\nu$ is a measure in $\mathcal {M}_{[0,1]^2}$
which satisfies $\int f(x) d\nu(x,y)=\int f(y) d\nu(x,y)$ for any  $f \in C[0,1]$.
 let $\nu=\theta P$ be the disintegration of  $\nu$.
To prove that $\nu$ belongs to $ \mathbf{M}$, we can use the fact that $\mathcal{A}$ is generated by the intervals, and thus  we just have to prove that $\theta(B) = \int P(x,B) d\theta(x) $
if $B$ is an interval. 

Therefore, Let  $B$ be an interval,  and $f_n\in C[0,1]$ a sequence of $[0,1]$-valued continuous functions that converges pointwise to $\chi_B$
(such a sequence always exists). By the dominated convergence theorem we have that
$$
\theta(B)=\int_{[0,1]} \chi_B(x) d\theta(x) = \lim_{n \rightarrow +\infty} \int_{[0,1]} f_n(x)d\theta(x)=$$
$$
=\lim_{n \rightarrow +\infty} \int_{[0,1]}  \int_{[0,1]} f_n(x) dP_x(y) d\theta(x)=
\lim_{n \rightarrow +\infty} \int_{[0,1]}  \int_{[0,1]} f_n(y) dP_x(y) d\theta(x)=$$
$$=\lim_{n \rightarrow +\infty} \int_{[0,1]} \varphi_n(x)d\theta(x)\,,$$
where $\varphi_n(x) \equiv \int_{[0,1]} f_n(y)dP_x(y)$. Now, defining $\varphi(x) \equiv \int_{[0,1]} \chi_B(y) dP_x(y)$, we can use again the dominated convergence theorem to get that $\varphi_n(x)  \rightarrow \varphi(x)$.
Hence the function  $\varphi_n$ is pointwise convergent and uniformly bounded. Using the dominated convergence theorem once more, we have that
$$
\theta(B)=\lim_{n \rightarrow +\infty} \int_{[0,1]} \varphi_n(x)d\theta(x) = \int_{[0,1]} \varphi(x) d\theta(x)= $$
$$=\int_{[0,1]}\int_{[0,1]}\chi_B(y)\; dP_x(y) d\theta(x)=
\int_{[0,1]}\int_{B} \;dP_x(y) d\theta(x)=$$
$$=  \int P(x,B) d\theta(x)   \,. $$

(b) 
Suppose $\nu_n \in  \mathbf{M}$, and $\nu_n \rightarrow \nu \in \mathcal {M}_{[0,1]^2}$ in the weak-$\star$ topology.
We have that $\int f d\nu_n \rightarrow \int f d\nu \; \forall f \in C([0,1]^2)$. In particular, if $f\in C[0,1]$, we have
$\int f(x) d\nu_n(x,y) \rightarrow \int f(x) d\nu(x,y)$ and $\int f(y) d\nu_n(x,y) \rightarrow \int f(y) d\nu(x,y)$.
Therefore, $\nu \in \mathbf{M}$ because
 $$\int f(x)-f(y) d\nu(x,y) = \lim_{n\rightarrow \infty} \int f(x)-f(y) d\nu_n(x,y) = 0\,,$$ for all $f \in C[0,1]$.

\cqd

The above formulation of the set $\mathbf{M}$ is more convenient for the duality of Fenchel-Rockafellar (see \cite{Roc} and the discussion on section \ref{section_maximprob}) required by proposition \ref{dualidade}.
It just says that both marginals in the $x$ and $y$ coordinates are the same.


 Sometimes we consider $\mu$ over
$[0,1]^\mathbb{N}$ and sometimes the corresponding projected  $\nu$ over $[0,1]^2$
(proposition \ref{relacao} below deals with projections of measures from  $ \mathcal M_0$ to $\mathbf{M}$).
We will
forget the word projected from now on, and the context will indicate which one we are working with.
 Note that, to make the lecture easier, we are using the following notation:
$\nu$ when we want to refer to a measure in $[0,1]^2$ and  $\mu$ for the measures in $\espaco$.

{\bf Remark :} We point out that maximizing $\int A d \nu$ for probabilities on $\nu\in \mathbf{M}$,
means a Kantorovich-Rubinstein (mass transhipment) problem where we assume the two marginals are the same (see \cite{Ra} Vol I section 4 for a related problem). The methods presented here can be used to get approximations of the optimal probability by absolutely continuous ones. These probabilities are obtained via the eigenfunctions of a Perron operator.

In the case we are analyzing, where the observable depends only on the two first coordinates, we will establish some connections between the measures in  $[0,1]^2$ and the measures in  $\espaco$,
 and we will see that the problem of maximization can be analyzed  as a problem of maximization among induced Markov measures in    $[0,1]^2$.

\begin{proposition}\label{relacao} Let   $A:\espaco\to
\re$ be a potential which depends only in the first two coordinates of  $\espaco$.  Then the following is true:

(a) There exists a map, not necessarily surjective, from
$\mathbf{M}$ to $ \mathcal M_0$.

\medskip

(b) There exists a map, not necessarily injective, from $
\mathcal M_0$ to $\mathbf{M}$.

\medskip

(c) $\displaystyle \max_{\mu\in\mathcal M_0 }\int A(x_1,x_2)\;
d\mu(\mathbf{x})=\max_{\nu\in  \mathbf{M}}\int A(x,y)\; d\nu(x,y)$

\end{proposition}

\textit{Proof: }
(a) A measure  $\nu \in \mathbf{M}$ can be disintegrated
as  $\nu=\theta P$, and then can be extended to a measure  $\mu\in \mathcal M_0$ by
\begin{equation}\label{st1} \mu(A_1...A_n):=\int_{A_1...A_n}\; dP_{x_{n-1}}( x_n)...dP_{x_1}( x_2) \;d\theta(x_1)  \, , \end{equation}

Also, we have
$$\int_{\espaco} A(x_1,x_2) d\mu (\mathbf{x})= \int_{[0,1]^2} A(x,y) d\nu(x,y) \,.$$

\medskip

(b) A measure   $\mu \in \mathcal {M}_{0}$ can be projected in a measure   $\nu \in \mathcal {M}_{[0,1]^2}$, defined by , for each Borel set   $B$ of $[0,1]^2$,
$$\nu(B) = \mu (\Pi^{-1}(B))\,,$$
where  $\Pi:\espaco \to [0,1]^2$ is the projection in the two first coordinates.
 Note that, by lemma  \ref{holon},
   $\nu \in \mathbf{M}$. Then we have  $$\int_{\espaco} A(x_1,x_2) d\mu (\mathbf{x})= \int_{[0,1]^2} A(x,y) d\nu(x,y) \,.$$

\medskip

(c) It follows easily by (a)  and (b). \cqd

\bigskip

\textbf{Remark:}   Note that in the item (a), in the particular case where   $\nu \in {\mathbf{M}_{ac}}$, we have that $\nu$
 can be disintegrated as $\nu=\theta K$, and then the stationary Markov  measure $\mu$ is given by
\begin{equation}\label{st2}  \mu(A_1...A_n):=\int_{A_1...A_n}K(x_{n-1},x_n)... K(x_1,x_2)\;\theta(x_1)\;dx_n...dx_1\, , \end{equation}
where  $A_1...A_n$ is a cylinder.


\section{The maximization problem}\label{section_maximprob}

We are interested in finding stationary Markov probabilities
$\mu_\infty$ on $ [0,1]^\mathbb{N}$ that maximize the value
$$ \int A(x_1,x_2) d \mu(\mathbf{x}),$$
over $\mathcal M_0$.

By  item  (c) of  proposition   \ref{relacao}: $\displaystyle \max_{\mu\in\mathcal M_0 }\int A\;
d\mu=\max_{\nu\in  \mathbf{M}}\int A\; d\nu$.

\bigskip

Hence, the  problem we are analyzing is equivalent   to the problem of  finding  $\nu_\infty$ which is  maximal for $\int A
d\nu$, among all $\nu \in    \mathbf{M}.$ Because once we have  $\nu_\infty$, by item
(a) of  proposition \ref{relacao}, we obtain a maximizing  Markov measure  $\mu_{\infty}$
among the holonomic measures.

\medskip

 As we only consider potentials of the form $A(x,y)$,  it is not
possible to have uniqueness of the maximizing measure on $\mathcal M_0$.  We just take into account  the
information of the measure on cylinders of size two. In any case,
the stationary Markov  probability we will describe below will also
solve this maximizing problem.

\bigskip

One of the main results we will get in this section is to be able to approximate singular probabilities by absolutely continuous probabilities (depending on a parameter $\beta$)  by means of eigenfunctions of a kind of Perron operator.

Now we will concentrate on the maximizing problem in $[0,1]^2$.

Let $A:[0,1] \times [0,1] \to \mathbb{R}$ be a continuous function. We will denote by
$$\mathfrak M_{0}:=\left\{\nu\in \mathbf{M} \;:\;\int A(x,y)\;d\nu(x,y)=m\right\}$$
where
$$m=\max_{\nu\in\ \mathbf{M}} \left\{ \int A(x,y) d\nu(x,y) \right\}\,.$$
A measure in $\mathfrak M_{0}$ will be called a maximizing measure on $ \mathbf{M}$.

\vspace{.5cm}

Consider now the variational problem

\begin{equation}\label{formterm}   \max_{\theta K\in\mathbf{M}_{ac}} \left\{ \int\beta  A(x,y)\theta(x) K(x,y) dx dy -
 \int \theta(x) K(x,y) \log\left(    K(x,y)      \right) dx dy \right\}
\end{equation}

In some sense we are considering above a kind of pressure problem (see \cite{PP}).

\begin{definition}
We define the {\em term of entropy} of an absolutely continuous
probability measure $\nu\in\mathcal M_{[0,1]^2}$, given by a density $\nu(x,y)dxdy$, as

\begin{equation}
S[\nu]=-\int \nu(x,y) \log\left(   \frac{\nu(x,y)}{\int\nu(x,z)dz}      \right) dx dy \,.
\end{equation}
\end{definition}

We remark that, in the case where  $A$ depends on all coordinates in $\espaco$,  the natural entropy (similar to Kolmogorov entropy for the case of the usual shift on the Bernoulli space) to be considered would be infinity. Therefore, it does not make sense to consider the associated concept of pressure (using Kolmogorov entropy)  and we believe it is not possible to go further in our reasoning to this more general setting. The bottom line is: we want to approximate singular probabilities by absolutely continuous probabilities (depending on a parameter $\beta$)  by means of eigenfunctions of a kind of Perron operator. We want to take limits in a parameter $\beta$ and this is easier to do  if we have a variational principle (like the one considered above).

It is easy to see that  any $\nu=\theta K\in \mathbf{M}_{ac}$  satisfies

\begin{equation}
S[\theta K]=-\int \theta(x) K(x,y) \log\left(    K(x,y)      \right) dx dy \,.
\end{equation}

 We call $S[\nu]= S[\theta K]$ the entropy penalized of the probability $\nu=\theta K\in \mathbf{M}_{ac}$.

\begin{lemma}\label{entropia}If $\nu=\theta K\in \mathbf{M}_{ac}$ and $K$ is   positive, then $S[\nu] \leq 0.$\end{lemma} 
\textit{Proof: } $\log$ is a concave function. Hence, by Jensen inequality, we have  $$-\int \theta(x) K(x,y) \log\left(    K(x,y)      \right) dx dy=\int \theta(x) K(x,y) \log\left( \frac{1}{K(x,y) }\right) dx dy\leq$$
$$\leq\log\int \theta(x) K(x,y)  \frac{1}{K(x,y) } dx dy=\log(1)=0.$$\cqd
\medskip

For each  $\beta$ fixed, we will exhibit  a measure $\nu_\beta$ in $\mathbf{M}_{ac}$ which maximizes  (\ref{formterm}). After, we will show that such $\nu_\beta$ will approximate in weak convergence   the probabilities $\nu_{\infty}$ which are maximizing for $A$ in the set $\mathbf{M}$.

In order to do that, we need to define the following operators:

\begin{definition} Let
${\cal L}_{\beta},{\bar{\cal L}}_{\beta}:C([0,1])\to C([0,1])$ be given by

\begin{equation}\label{L}   {\cal L}_{\beta} \varphi (y) =\int e^{\beta A(x,y)} \, \varphi(x)dx,  \end{equation}

\begin{equation}\label{bar L}   {\bar{\cal L}}_{\beta} \varphi (x) =\int e^{\beta A(x,y)} \, \varphi(y)dy.  \end{equation}
\end{definition}
\vspace{0.2cm}

We refer the reader to \cite{Ka} and \cite{Sch} chapter IV for a general reference on positive integral operators.
\vspace{0.2cm}

The above definitions are quite natural and extend the usual Ruelle-Perron operator definition. In the present situation the state space is continuous and an integral should take place of the sum. We are interested in approximating singular measures (which are maximizing for $A$) by absolutely continuous probabilities, therefore, it is natural to integrate with respect to Lebesgue measure.

\begin{theorem}\label{KreinRuthman} 
The operators  ${\cal L}_{\beta}$ and ${\bar{\cal L}}_{\beta}$ have the same  positive maximal eigenvalue $\lambda_{\beta}$, which is simple and isolated. The eigenfunctions associated are  positive functions.
\end{theorem}

\textit{Proof:}
We can see that
${\cal L}_{\beta}$ is a compact operator, because the image by
${\cal L}_{\beta}$ of the unity closed ball of $C([0,1])$ is a
equicontinuous family in $C([0,1])$: we know that $e^{\beta A}$ is a
uniformly continuous function, and then, if $\varphi$ is in the
closed unit ball, we have
$$
|{\cal L}_{\beta} \varphi (y) - {\cal L}_{\beta} \varphi (z)| \leq \int |e^{\beta A(x,y)}-e^{\beta A(x,z)}| \,|\varphi(x)| dx  \leq |e^{\beta A(x,y)}-e^{\beta A(x,z)}| < \delta\, ,
$$
if, $y$ and $z$ are close enough. Thus, we can use Arzela-Ascoli Theorem to prove the compactness of ${\cal L}_{\beta}$ (see also \cite{Sch} Chapter IV section 1).

The spectrum of a compact operator is a sequence of eigenvalues that converges to zero, possibly added by zero. This implies that any non-zero eigenvalue of ${\cal L}_{\beta}$ is isolated
(i.e. there are no sequence in the spectrum of ${\cal L}_{\beta}$ that converges to some non-zero eigenvalue).

The definition of ${\cal L}_{\beta}$ now shows that it preserves the cone of positive functions in $C([0,1])$, indeed, sending a point in this cone to the interior of the cone. This means that ${\cal L}_{\beta}$ is a  positive operator.

The Krein-Ruthman Theorem (Theorem 19.3 in \cite{De}) implies that there exists a
 positive eigenvalue $\lambda_{\beta}$, which is  maximal (i.e. $\lambda_{\beta}>|\lambda|$, if $\lambda\neq\lambda_{\beta}$ is in the spectrum of ${\cal L}_{\beta}$)
and simple  (i.e. the eigenspace associated to $\lambda_{\beta}$ is one-dimensional). Moreover $\lambda_{\beta}$ is associated to a positive eigenfunction
$\varphi_{\beta}$. 

If we proceed in the same way, we get the same conclusions about the operator ${\bar{\cal L}}_{\beta}$, and we get
the respective eigenvalue $\bar\lambda_{\beta}$ and eigenfunction $\bar{\varphi}_{\beta}$.

In order to prove that $\bar\lambda_{\beta}= \lambda_{\beta}$, we use
the positivity of $\varphi_{\beta}$ and $ \bar{\varphi}_{\beta}$ and the fact that ${\bar{\cal L}}_{\beta}$
is the adjoint of ${{\cal L}}_{\beta}$ (here we see that  our operators can be, in fact, defined in the Hilbert space $L^2([0,1])$, which contains $C([0,1])$). We have
$<\varphi_{\beta},\bar{\varphi}_{\beta}> = \int \varphi_{\beta}(x) \bar{\varphi}_{\beta}(x) dx > 0$, and
$$
\lambda_{\beta} <\varphi_{\beta},\bar{\varphi}_{\beta}> = <{{\cal L}}_{\beta}\varphi_{\beta},\bar{\varphi}_{\beta}> =
<\varphi_{\beta},{\bar{\cal L}}_{\beta}\bar{\varphi}_{\beta}> = \bar\lambda_{\beta} <\varphi_{\beta},\bar{\varphi}_{\beta}>.$$ \cqd

\medskip

An estimate on the spectral gap for the operator ${\cal L}_\beta$, where $\beta>0$, is given in \cite{Os} \cite{Hop}: suppose
$$ \tilde{M}= \sup_{(x,y)\in [0,1]^2}\, A(x,y),\,\,\,\,\,\,\,\text{and}\,\,\,\,\,\tilde{m}= \inf_{A(x,y)\in [0,1]^2}\, A(x,y).$$

If $\lambda_\beta$ is the main eigenvalue, then, by theorem 4 of \cite{Hop}, any other $\lambda$ in the spectrum of
${\cal L}_ \beta$ satisfies
$$ \lambda_\beta  \,\,(\frac{\tilde{M}^\beta - \tilde{m}^\beta}{\tilde{M}^\beta + \tilde{m}^\beta})\,\, >\,\,  \lambda.$$

With this information one can give an estimate of the decay of correlation for functions evolving under the probability of the Markov Chain associated to such value $\beta$ (see next proposition). The proof of this claim is similar to the reasoning in chapter 2 page 26 in \cite{PP}, which deals with the case where the state space is discrete.

\bigskip

Let us call $\varphi_{\beta},\bar{\varphi}_{\beta}$  the
positive eigenfunctions for ${\cal L}_{\beta}$ and ${\bar{\cal
L}}_{\beta}$ associated to $\lambda_{\beta}$, which satisfy $\int
\varphi_{\beta}(x)\,dx=1$ and $\int \bar\varphi_{\beta}(x)\,dx=1$.

\vspace{0.3cm}

We will define a density  $\theta_{\beta}:[0,1]\to \re$ by \begin{equation}\label{theta}\theta_{\beta}(x):=\frac{\varphi_{\beta}(x)\,\, \bar\varphi_{\beta}(x)}{\pi_{\beta}},\end{equation} where $\pi_{\beta}=\int \varphi_{\beta}(x) \bar\varphi_{\beta}(x)dx$,
and a transition $ K_{\beta}:[0,1]^2\to \re$ by
\begin{equation}\label{K}K_{\beta}(x,y):=\frac{e^{\beta A(x,y)}\,\,\bar \varphi_{\beta}(y)}{\bar \varphi_{\beta}(x)\,\lambda_{\beta}}. \;\end{equation}

Let $\nu_{\beta}\in \mathcal {M}_{[0,1]^2}$ be defined by
\begin{equation}\label{nubeta}d\nu_{\beta}(x,y):=\theta_{\beta}(x)
K_{\beta}(x,y)dxdy.\end{equation}

It is easy to see  that $ \theta_{\beta},
K_{\beta}$ satisfy equations \eqref{transicao}, \eqref{normalizada} and \eqref{equilibrio},  hence $\nu_{\beta}\in \mathbf{M}_{ac}$.

\begin{proposition}\label{maxMarkov} The Markov measure $\nu_{\beta}=\theta_{\beta}K_{\beta}$ defined above maximize
$$  \int   \, \beta\,  A(x,y)\;\theta(x) K(x,y) dx dy -
 \int \theta(x) K(x,y) \log\left(    K(x,y)      \right) dx dy$$
 over all absolutely continuous  Markov measures. Also

 $$\log \lambda_{\beta}=\int \beta A \;\theta_{\beta} K_{\beta}dxdy +S[\theta_{\beta} K_{\beta}].$$
\end{proposition}

\textit{Proof: }
 By the definition of the functions $\theta_{\beta}, K_{\beta}$, we have
$$S[\theta_{\beta} K_{\beta}]=-\int (\beta A(x,y)+ \log {\bar \varphi_{\beta}(y)}-\log{\bar \varphi_{\beta}(x)-\log \lambda_{\beta}}  ) d\nu_{\beta}\,.  $$

Then
$$\int \beta A \;\theta_{\beta} K_{\beta}dxdy +S[\theta_{\beta} K_{\beta}]=$$
$$=\log \lambda_{\beta}+\int (\log\bar\varphi_{\beta}(x)-\log{\bar\varphi_{\beta}(y)})\theta_{\beta}(x)
K_{\beta}(x,y)dxdy\,,$$ 
and the last integral is zero because $\nu_{\beta}=\theta_{\beta}K_{\beta}  \in \mathbf{M}_{ac}$.




To show that $\nu_{\beta}$ is maximizing let $\nu$ be any measure in $\mathbf{M}_{ac}$ and $0\leq\tau\leq 1$. We claim that the function

$$I[\tau]:= \int {\beta}A d\nu_{\tau}+ S[\nu_{\tau}]$$ where $\nu_{\tau}=(1-\tau)\nu_{\beta}+\tau \nu$, is concave and $I'(0)=0$

Indeed, see proof of theorem 33 of \cite{GV}. We just  point out  that the entropy term in \cite{GV} has a difference of sign.\cqd

\bigskip

\begin{lemma}\label{relcompact} (a) There exists a constant $c>0$ such that, for all $x\in[0,1]$, we have

$$e^{-\beta c}\leq \varphi_{\beta}(x)\leq e^{\beta c} \;\;\;\mbox{and}\;\;\;\; e^{-\beta c}\leq \bar\varphi_{\beta}(x)\leq e^{\beta c}.$$
Also,
$$ \beta \mapsto \frac{1}{\beta}\log\pi_{\beta} \;\;\; \mbox{and} \;\; \;\; \beta \mapsto  \frac{1}{\beta}\log\lambda_{\beta}  $$
are bounded functions, defined for $\beta>0$.

(b) The sets $$\left\{ \frac{1}{\beta}\log(\varphi_{\beta}) \;\;| \; \beta>1 \, \right\} \;\;\mbox{and} \;\;
\left\{ \frac{1}{\beta}\log(\bar\varphi_{\beta}) \;\;| \; \beta>1 \,
\right\}$$ are equicontinuous, and relatively compact in the supremum norm. 


\end{lemma}

 \textit{Proof: }(a)
 Fix $\beta>0$. Using the normalization $\int \varphi_{\beta}(z)dz=1$, we choose $x_0$ and $x_1$ in $[0,1]$ satisfying  $\varphi_{\beta}(x_0)\leq 1$ and $\varphi_{\beta}(x_1) \geq 1$.
  Now, if $\|A\|$ is the supremum norm of $A$, we have
 $$\lambda_{\beta}=\frac{1}{\varphi_{\beta}(x_1)} \int e^{\beta A(z,x_1)} \varphi_{\beta}(z)dz \leq e^{\beta \|A\|} \;\mbox{and}$$
 $$\lambda_{\beta}=\frac{1}{\varphi_{\beta}(x_0)} \int e^{\beta A(z,x_0)} \varphi_{\beta}(z)dz \geq e^{-\beta \|A\|} \; .$$
 Thus, $-\|A\| < \frac{1}{\beta}\log\lambda_{\beta} < \|A\|$.

\medskip
Now we use the inequalities above and the fact that
 $$\varphi_{\beta}(x)=\frac{1}{\lambda_{\beta}} \int e^{\beta A(z,x)} \varphi_{\beta}(z)dz $$
to prove that 
 $$\varphi_{\beta}(x) \leq
 \frac{1}{\lambda_{\beta}} \int e^{\beta \|A\|} \varphi_{\beta}(z)dz \leq e^{2 \beta \|A\|}\,,$$
and
$$\varphi_{\beta}(x) \geq
 \frac{1}{\lambda_{\beta}} \int e^{-\beta \|A\|} \varphi_{\beta}(z)dz  \geq e^{-2\beta \|A\|}\,.$$
We define $c= 2 \|A\|$. The eigenfunctions $\overline{\varphi}_{\beta}$ are bounded by an analogous estimative.
Now, $\pi_{\beta}= \int \varphi_{\beta}(x) \overline{\varphi}_{\beta}(x) dx$, and thus $  e^{-2\beta c} \leq \pi_{\beta} \leq e^{2\beta c}$, which implies that $\frac{1}{\beta}\log\pi_{\beta}$ is a bounded function of $\beta$.

\medskip

(b) We just have to prove the equicontinuity of both sets. Once we have that, and considering
the fact that both sets are sets of functions defined in the compact set $[0,1]$, we use item (a) and Arzela-Ascoli's Theorem to get the relative compactness of these sets.


To have the equicontinuity for the first set, let $y$ be a point in $[0,1]$, and let $\beta >1$. Let $\epsilon >0$.
We will use the fact that $A$ is a uniformly continuous map: We know there exists $\delta>0$, such that $|y-z|<\delta$, implies
$|A(x,y)-A(x,z)|<\epsilon, \; \forall x\in[0,1]$.
Without any loss of generality, we suppose that $\varphi_{\beta}(y) \geq \varphi_{\beta}(z)$.
We have:
$$ \left|\frac{1}{\beta}\log(\varphi_{\beta}(y)) -
\frac{1}{\beta}\log(\varphi_{\beta}(z))  \right|
= $$
$$ = \frac{1}{\beta}
\left( \log\left(\frac{1}{\lambda_{\beta}}
\int e^{\beta A(x,y)} \varphi_{\beta}(x) dx \right)
- \log\left(\frac{1}{\lambda_{\beta}}
\int e^{\beta A(x,z)} \varphi_{\beta}(x) dx \right)
\right) =
$$
$$ = \frac{1}{\beta}
 \log\left(
\frac
{\int e^{\beta A(x,y)} \varphi_{\beta}(x) dx }
{\int e^{\beta A(x,z)} \varphi_{\beta}(x) dx }\right)
\leq
\frac{1}{\beta}
 \log\left(
\frac
{\int e^{\beta (A(x,z)+\epsilon )} \varphi_{\beta}(x) dx }
{\int e^{\beta A(x,z)} \varphi_{\beta}(x) dx }\right) =
$$
$$ =
\frac{1}{\beta}
 \log\left(
e^{\beta \epsilon}
\frac
{\int e^{\beta A(x,z)} \varphi_{\beta}(x) dx }
{\int e^{\beta A(x,z)} \varphi_{\beta}(x) dx }\right)
=
\epsilon\,.$$

We  prove the equicontinuity for the second set in the same way.

\cqd

From the above, we can find $\beta_n\to \infty$ which defines convergent subsequences $\frac{1}{\beta_n}\log\varphi_{\beta_n}.$

Let us fix a subsequence $\beta_n$ such that $\beta_n\to \infty$ and  all the three following limits exist:
$$V(x):=\lim_{n\to \infty}\frac{1}{\beta_n}\log\varphi_{\beta_n}(x)\;\;\;, \;\;\;\bar V(x):=\lim_{n\to \infty}\frac{1}{\beta_n}\log\bar\varphi_{\beta_n}(x)$$
$$\tilde m:=\lim_{n\to \infty}\frac{1}{\beta_n}\log \lambda_{\beta_n}$$
Note that the limits defining $V$ and $\bar V$ converge uniformly. In principle, the function $V$ depends on the sequence $\beta_n$  we choose.

\begin{proposition}[Laplace's Method]Let $f_{k}:[0,1]\to\re$ be a sequence of functions that converge uniformly, as $k$ goes to $\infty$, to a function $f:[0,1]\to\re$.
Then
$$\lim_{k}\frac{1}{k}\log \int_0^1 e^{k f_{k}(x)} dx = \sup_{x\in [0,1]} f(x)$$

\end{proposition}

\begin{lemma}
$$\lim_{n\to \infty}\frac{1}{\beta_n}\log \pi_{\beta_n}=\max_{x\in [0,1]}(V(x)+\bar V(x))$$
\end{lemma}
  \textit{Proof: }
  $$\pi_{\beta_n}=\int_0^1\varphi_{\beta_n}(x) \bar\varphi_{\beta_n}(x)dx=\int_0^1 e^{\beta_n(\frac{1}{\beta}\log\varphi_{\beta_n}(x)+\frac{1}{\beta} \log\bar\varphi_{\beta_n}(x))}dx$$
And note that $\frac{1}{\beta_n}\log\varphi_{\beta_n}(x)\to V(x), \frac{1}{\beta_n}\log\bar \varphi_{\beta_n}(x)\to \bar V(x)$ uniformly. Hence it follows by Laplace's Method.\cqd

\vspace{.5cm}

Also by Laplace's method we have the following lemma:

\begin{lemma}\label{subaction2}$$V(y)=\max_{x\in [0,1]}(V(x)+A(x,y)-\tilde m)$$
and $$\bar V(x)=\max_{y\in [0,1]}(\bar V(y)+A(x,y)-\tilde m).$$
\end{lemma}
\vspace{.5cm}

For some subsequence (of the subsequence $\{\beta_n\}$ fixed after the proof of lemma \ref{relcompact}, which we will also denote  by $\{\beta_n\}$), the measures  $\nu_{\beta_n}$ defined in (\ref{nubeta}) weakly converge to a measure  $\nu_{\infty}\in \mathcal {M}_{[0,1]^2}$. Then
$$\lim_{n\to \infty}\int A   d\nu_{\beta_n}= \int A   d\nu_{\infty}\,.$$

\begin{lemma}\label{markov}  The measure $\nu_{\infty}\in \mathbf{M}$.\end{lemma} \textit{Proof: }
 As $\nu_{\beta_n}\in\mathbf{M}_{ac}\subset\mathbf{M} $, by item (b) of lemma \ref{holon} we have that $\nu_{\infty}\in \mathbf{M}$.\cqd

\begin{theorem} \label{max} 
$$\int A(x,y)d\nu_{\infty}(x,y)=m$$ i.e.,  $\nu_{\infty}$ is a maximizing measure on $ \mathbf{M}$.
\end{theorem}

In order to prove theorem \ref{max} we need first some new results.

\begin{proposition}\label{dualidade} Given a   potential $A\in C([0,1]^2)$ ,  we have that
$$\sup_{\nu\in\ \mathbf{M}}\int A d\nu=\inf_{f\in\; C([0,1])}\;\max_{(x,y)}\left(A(x,y)+f(x)-f(y)\right)$$
\end{proposition}

This proposition will be a consequence of the Fenchel-Rockafellar
duality theorem (see \cite{Roc}). Let us fix the setting we consider in order to  apply this theorem.

Let $C([0,1]^2)$ be the set of continuous functions in $[0,1]^2$ with the supremum norm and
$\mathcal S$ the set of signed measures over the Borel
$\sigma-$algebra of $[0,1]^2$.

Consider the convex correspondence $H:C([0,1]^2)\to\re$ given by $H(\phi)=\max(A+\phi)$ and
$$\mathcal C:=\{ \phi\in C([0,1]^2) \;:\; \phi(x,y)=f(x)-f(y) \mbox {\;\;, for some\;\; }f\in C([0,1])\}$$

We define a concave correspondence $G:C([0,1]^2)\to\re\cup
\{-\infty\}$ by $G(\phi)=0$ if $\phi\in\mathcal C$ and
$G(\phi)=-\infty$ otherwise.

Then the corresponding Fenchel transforms,

$H^*:\mathcal S\to \re \cup \{+\infty\}$, $G^*:\mathcal S\to \re
\cup \{-\infty\}$, are given by
$$H^*(\nu)=\sup_{\phi\in C([0,1]^2)}\left[ \int \phi(x,y)d \nu(x,y)
-H(\phi)\right]$$

and
$$G^*(\nu)=\inf_{\phi\in C([0,1]^2)}\left[ \int \phi(x,y)d \nu(x,y)
-G(\phi)\right]$$

We define $\mathcal S_0:=\{\nu\in\mathcal S \;:\; \int f(x)-f(y)d\nu(x,y)=0 \;\;\; \forall f\in C[0,1]\} $, and we note that $\mathcal S_0\cap \mathcal M_{[0,1]^2}={\mathbf{M}} $.

\begin{lemma}\label{legendre}Given $H$ and $G$ as above, then

$$
 H^*(\nu)= \left\{
\begin{array}{ccc}
-\int A(x,y)d\nu(x,y)&\mbox{if}&  \nu\in \mathcal M_{[0,1]^2}\\
+\infty& & \mbox{otherwise}
\end{array}
\right.
$$

$$
 G^*(\nu)= \left\{
\begin{array}{ccc}
0 &\mbox{if}&  \nu\in  \mathcal S_0\\
-\infty& & \mbox{otherwise}
\end{array}
\right.
$$

\end{lemma}
This lemma follows from lemma 2 of \cite{GL}.

\vspace{.4cm}
 \textit{Proof of proposition \ref{dualidade}:}
The duality theorem of Fenchel-Rockafellar says that

$$\sup_{\phi\in C([0,1]^2)}[G(\phi)-H(\phi)]=\inf_{\nu\in\mathcal
S}[H^*(\nu)-G^*(\nu)]\,.$$

Hence, by lemma \ref{legendre} and the uniform convergence we have

$$\sup_{\phi\in\;\mathcal C}[-\max_{(x,y)}(A+\phi)(x,y)]=\inf_{\nu\in \mathbf{M}}[-\int A(x,y)
d\nu(x,y)]\,.$$ Using the definition of $\mathcal C$  we have that

$$\sup_{\nu\in \mathbf{M}}\int A d\nu=\inf_{f\in\; C([0,1])}\;\max_{(x,y)}
\left(A(x,y)+f(x)-f(y)\right)\,.$$ \cqd

 \begin{lemma}
\label{migualmtilde} $\tilde m=m$.
\end{lemma}

 \textit{Proof: }
 Note that by proposition \ref{dualidade} and lemma \ref{subaction2} we have that $m\leq \tilde m$.
To show the other inequality remember that
$$\log \lambda_{\beta_n}=\int \beta_n A \;d\nu_{\beta_n}+S[\nu_{\beta_n}].$$
Then
$$\tilde m=\lim_{n\to \infty}\int  A \;d\nu_{\beta_n}+\frac{1}{\beta_n}S[\nu_{\beta_n}].$$
Note that  $\nu_{\beta_n}\in  \mathbf{M}$, which implies  $\int A \;d\nu_{\beta_n}\leq m $.

As  $S[\nu_{\beta_n}]\leq 0$,   we have
$$\int  A \;d\nu_{\beta_n}+\frac{1}{\beta_n}S[\nu_{\beta_n}]\leq m\;\;\;\;\; \forall n$$
Then
$\tilde m \leq m$.
\cqd

\vspace{.7cm}
\textit{Proof of Theorem \ref{max}:}
Remember that $\nu_{\beta_n}\rightharpoonup \nu_{\infty}$, then

$$\lim_{n\to\infty}\int A d\nu_{\beta_n} = \int A d\nu_{\infty}.$$  By lemma \ref{migualmtilde} and the fact that $S[\nu_{\beta_n}]\leq 0$, we obtain
$$ m=\lim_{n\to \infty}\int  A \;d\nu_{\beta_n}+\frac{1}{\beta_n}S[\nu_{\beta_n}]\leq \int  A \;d\nu_{\infty}.$$ Hence, using lemma \ref{markov}, we have that
$m=\int A d\nu_{\infty}.$\cqd

\section{Uniqueness of maximizing measures and calibrated subactions}\label{uniqueness}

We want to remark here  that for the results of  this section we
were inspired
 by the  works of \cite{Gom}, \cite{GLM} and
 \cite{GL}. Hence, jointing all these ideas, and what was proved before, we are able to show that
 there is a unique maximizing probability for $A$ in $\mathbf{M}$, if $A$ is  generic in Ma\~n\'e's sense, the potential $A$ is $C^2$
 and satisfy the twist property. Similar result is true for the calibrated subaction.
The precise definitions will be given in what follows.

The differentiable structure of $[0,1]$ will help us  to get the uniqueness required when we want to show   the graph property for the support of the maximizing probability.

We repeat here the important definition (Definition \ref{sub}) of
forward (backward)-calibrated subactions:

\begin{definition}
A continuous function $u: [0,1] \to \mathbb{R} $ is called a

(a) {\em  calibrated forward-subaction} if, for any $y$ we have
\begin{equation}\label{c} u(y)=\max_x [A(x,y)+ u(x)-m].\end{equation}

(b) {\em calibrated backward-subaction} if, for any $x$ we have
\begin{equation}\label{d} u(x)=\max_y [A(x,y)+ u(y)-m].\end{equation}
\end{definition}

Note also that if we add a constant to a calibrated forward-subaction, this will be a new calibrated forward-subaction. When we say here that under some conditions, the calibrated  forward-subaction is unique, we say this up to an additive constant.

Note that $V$ and $\bar V$ defined in lemma \ref{subaction2} are,
respectively, forward and backward calibrated subactions
  (remember that $\tilde m =m $ by lemma \ref{migualmtilde}).

Subactions (see also \cite{CLT}) play the role in discrete time dynamics of fixed points of the Lax-Oleinik operators of Mather Theory \cite{Fathi}.

\bigskip

Let $u$ be a calibrated backward-subaction, using the fact that $[0,1]$ is compact, there exists  $y(x)$ (maybe not unique) such that

\begin{equation}\label{b}u(x)= A(x,y(x))+u(y(x))-m.\end{equation}

\begin{proposition}\label{cohomologia}
  Let $\nu\in\mathfrak M_{0}$ be any maximizing measure, and $u$ be a calibrated backward-subaction. Then for all $(x,y)\in\supp (\nu)$ we have
$$ u(x)= A(x,y)+u(y)-m .$$
\end{proposition}

\textit{Proof:} Note that  $ u(x)\geq A(x,y)+u(y)-m $ for all $y\in[0,1]$.  As $\nu\in\mathfrak M_{0}$, we have $\int A d\nu=m$ and  $\int (u(x)-u(y))\,d\nu=0$. This proves that the equality in the statement of the theorem is true $\nu$-almost everywhere, in the other points of the support of $\nu$ this holds by continuity.
\cqd
\vspace{0.2cm}

We point out that a calibrated-subaction (backward or forward) does  not need  to be differentiable. We want to show that, in certain  points of $[0,1]$, a calibrated-subaction is  differentiable. In order to do that we introduce the following generalized differentials.

\begin{definition}  Let $u:B\to\re$ and $x\in B$, where $B$ is an open set in $[0,1]$. The sets
$$D^+u(x)=\{\,p \in \mathbb{R} \,| \, \limsup_{|v|\to 0} \, \frac{u(x+v) - u(x) - p \,v}{|v|} \leq 0\,\},$$

$$D^-u(x)=\{\,p \in \mathbb{R} \,| \, \liminf_{|v|\to 0} \, \frac{u(x+v) - u(x) - p \,v}{|v|} \geq 0\,\},$$
are called, respectively, the {\em superdifferential} and the {\em subdifferential of $u$ at $x$}.

\end{definition}

The main point here is that the differentiable structure of $[0,1]$ will help us  to get the uniqueness required by what we will call later the graph property.

\begin{proposition}\label{super} Let $u:B\to\re$ and $x\in B$, where $B$ is an open set in $[0,1]$.  $D^-u(x)$ and $D^+u(x)$ are both nonempty sets if and only if $u$ is differentiable at $x$. In this case $D^-u(x)=D^+u(x)=Du(x)$.

\end{proposition}
\textit{Proof:} See proposition 3.1.5 of \cite{CS}.\cqd

\vspace{0.2cm}

\begin{lemma} \label{subdif} Let  $u$ be a  calibrated backward-subaction.  We have the following statements:
 $$(a)\;\;\;\forall x\in[0,1], \;\;\; D^-u(x)\neq \emptyset\;\; ,\;\;\;\; \mbox{ and } \;\;\;\;\;\frac{\partial A}{\partial x}(x,y(x)) \in D^-u(x),$$ where $ y(x)$ is such  that $(x,y(x))$ satisfies  equation (\ref{b});

\bigskip
For $(x,y(x))$ satisfying  equation (\ref{b}):

\bigskip

$ (b)\;\;\;  \displaystyle{D^+u (y(x))\neq \emptyset \;\;\;\; \mbox{ and } \;\;\;-\frac{\partial A}{\partial y}(x,y(x))\in D^+u(y(x)); }$

\;\;\;\; (c)\;\; $u$ is differentiable at $y(x)$.

\end{lemma}

\textit{Proof:} (a) Let $x\in [0,1]$, then there exists $ y(x)$ such that  $(x,y(x))$ satisfies  equation (\ref{b}). For any $w\in [0,1]$, using equation \eqref{d}, we have that
$$u(x+w)\geq A(x+w,y(x))+u(y(x))-m\,.  $$
This inequality and the equality in equation \eqref{b}, give  that

$$u(x+w)-u(x)- A(x+w,y(x))+A(x,y(x)) \geq 0 ,$$ and then
$$\liminf_{|w|\to 0}\frac{u(x+w)-u(x)- ( \frac{\partial A}{\partial x}(x,y(x))\;w+o(w))}{|w|} \geq 0 \,,$$ and this implies $\frac{\partial A}{\partial x}(x,y(x)) \in D^-u(x)$.

(b) Also for $(x,y(x))$ satisfying  (\ref{b}) and any $w\in [0,1]$, using equation \eqref{d}, we have
$$ u(x)\geq A(x,y(x)+w)+u(y(x)+w)-m.$$ Using equation (\ref{b}), we get that
$$u(y(x)+w)-u(y(x))+A(x,y(x)+w)-A(x,y(x))\leq 0. $$ Then
$$\limsup_{|w|\to 0}\frac{u(y(x)+w)-u(y(x))-(-\frac{\partial A}{\partial y}(x,y(x))\; w+o(w))}{|w|}\leq 0.$$
Hence, $-\frac{\partial A}{\partial y}(x,y(x))\in D^+u(y(x))$.

(c) It is just a consequence of items (a) and (b) and of proposition \ref{super}.
\cqd


\begin{lemma}\label{z} For any measure    $\nu \in \mathbf{M}$, we have that, for almost every point $(x,y)\in\supp(\nu)$, there exists $z$ such that $(z,x)\in\supp(\nu)$.

\end{lemma}
\textit{Proof:}  Define the set
$$R=\{(x,y)\in \supp(\nu)\; : \; \supp(\nu)\cap([0,1]\times\{x\})=\emptyset \}$$
Suppose, by contradiction, that  $\nu(R)=\epsilon >0$.

Let $\pi_j: [0,1]^2 \to [0,1]$ be the projection on the $j$-th coordinate.

Let $\nu_2$ be the measure on the Borel sets of $[0,1]$ given by $\nu_2(B) = \nu (\pi_2^{-1}(B))$, where $B$ is any Borel set in $[0,1]$.

Consider $R_1=\pi_1(R)$. We have
$$R_1=\{x\in \pi_1(\supp(\nu))\; : \; \supp(\nu)\cap([0,1]\times\{x\})=\emptyset\}.$$
We claim that
$$\nu_2(R_1) =
\int_{\supp (\nu)}\chi_{R_1}(y)d\nu_{\infty}(x,y)=0.$$
Indeed, the first equality is immediate. To prove the second equality,
take  $(x,y)\in\supp(\nu)$. We have two possibilities:
If $y\notin\pi_1(\supp(\nu))$, then $y\notin R_1$.
And if  $y\in\pi_1(\supp(\nu))$ we have $(x,y)\in\supp(\nu)$ and then $y\notin R_1$. This shows the claim.

By the other hand, note that $\pi_1^{-1}(R_1)\cap \supp(\nu)=R$, and thus
$$\int_{\supp (\nu)}\chi_{R_1}(x)d\nu(x,y)=\int_{\supp (\nu)}\chi_{\pi_1^{-1}(R_1)}(x,y)d\nu(x,y)=\nu(R)=\epsilon.$$

Now let $U$ be an open set of $[0,1]$ which contains $R_1$ and such that $\nu_2(U) < \nu_2(R_1) + \epsilon/2 =
\epsilon/2$.
Consider a sequence of continuous  function $f_j$ such that $f_j \uparrow \chi _{U}$.
Using the monotonous convergence theorem and $\nu\in \mathbf{M}$, we have:

$$\epsilon/2 > \nu_2(U) = \int\chi_{U}(y)d\nu(x,y) = \lim_j \int f_j(y)d\nu(x,y)=$$
$$= \lim_j\int f_j(x)d\nu(x,y)=\int \chi_{U}(x)d\nu(x,y) \geq
\int \chi_{R_1}(x)d\nu(x,y)=\epsilon$$ which is a contradiction.\cqd

\begin{theorem}\label{grafico}Let $\nu\in\mathfrak M_{0}$ be any maximizing measure. If the observable $A$ is $C^2$, and $\frac{\partial^2 A}{\partial x\partial y}>0$, then  the measure
 $\nu$ is supported on a graph.
\end{theorem}

\textit{Proof:} Let $u$ be any
 calibrated backward-subaction and  $(x_0,y_0)\in\supp (\nu)$,  then by proposition \ref{cohomologia} $(x_0,y_0)$ satisfies equation (\ref{b}).

 On the other hand,  by  lemma \ref{z},
  there exists $z_0$ such that $(z_0,x_0)\in\supp (\nu)$, and this  means that $x_0=y(z_0)$. Thus  item (c) of  lemma \ref{subdif} implies that $u$ is differentiable at $x_0$.

   Now by item (a) of lemma \ref{subdif} and proposition \ref{super}, we have that  \begin{equation}\label{dif}Du(x_0)=\frac{\partial A}{\partial x}(x_0,y_0)\end{equation}

   Note that, for any fixed $p$ and $x$,  the equation
$p=\frac{\partial A}{\partial x}(x,y)$ has at most  one solution
$y(x,p)$  because  $\frac{\partial^2 A}{\partial x\partial y}>0$.
  Let $p=Du(x_0)$ and $x=x_0$ then  $y_0$ is the unique point that satisfies the equation \eqref{dif}.  \cqd

  \bigskip

 \textbf{Remark:} Using the same arguments of the   proof of theorem \ref{grafico}, we see that,   if $u$ is a calibrated subaction, $u$ is differentiable at $x$ and  $(x,y)$ satisfies equation (\ref{b}) (note that, for each $x$ there exists at least one $y$ with this property), then we  get that this $y$ is the unique point that satisfies the equation $Du(x)=\frac{\partial A}{\partial x}(x,y)$. Therefore, $y$ is  the unique point that satisfies equation (\ref{b}).
\vspace{.3cm}

\begin{lemma}\label{uniaosuportes}If the observable $A$ is $C^2$, and $\frac{\partial^2 A}{\partial x\partial y}>0$, then $\cup_{\nu \in \mathfrak M_{0}}  \supp(\nu)$  is contained in a graph.
\end{lemma}
\textit{Proof:}
Let $\nu_1$ and $\nu_2$ be two maximizing measures. Suppose there exists $x \in \pi_1(\supp(\nu_1)) \cap \pi_1(\supp(\nu_2))$.
Let $y_1$ and $y_2$ be the (unique) points such that $(x,y_1) \in \supp(\nu_1)$ and  $(x,y_2) \in \supp(\nu_2)$.

Let $u$ be a calibrated backward-subaction, using the same arguments of the proof of theorem \ref{grafico} for $(x,y_1) \in \supp(\nu_1)$, and for $(x,y_2) \in \supp(\nu_2)$, we get, respectively
$$Du(x_0)=\frac{\partial A}{\partial x}(x_0,y_1) \mbox{\;\;\;\;and\;\;\;}Du(x_0)=\frac{\partial A}{\partial x}(x_0,y_2)\,.$$

But, as before, the equation  $p=\frac{\partial A}{\partial x}(x,y)$ has at most one solution $y(x,p)$, then $y_1=y_2$. \cqd

\begin{definition}
 Given $k$ and  $ x ,
y \in [0,1]  $, we will call a {\em $k$-path beginning in $x $ and ending
at $ y$} an ordered sequence of points
$$ (x_1,....,x_k)\in [0,1]\times ...\times [0,1] $$
satisfying $ x_1 = {x} $,  $x_k=y$.
\end{definition}

We will denote by $ \mathcal P_k (x, y)$ the set of such $k$-paths.
\hspace{.1in}

\textbf{Remark:} 1) Here we shall note that the results that we will get  can not be a particular case of the results obtain in \cite{Gom}, \cite{GLM} for  the theory of Aubry-Mather, because in A-M theory  a Lagrangian $L:[0,1]\times\re\to\re$, satisfy the hypothesis  that  $L(x,v)\to +\infty$ when $|v|\to\infty$.

2) A path in A-M theory (see \cite{GLM}) is an orderer  sequence of points  $ (x_0,....,x_k)\in \rn\times ...\times \rn $ such that for each  $x_j$ we associate a velocity $v_j=x_{j+1}-x_j$, $0\leq j<k$. With  those pairs  $(x_j,v_j)$ we are able to calculate the action of the path $ (x_0,....,x_k)$. In our setting there is no velocity and only the points of the path are used to calculate the action of the path.

\bigskip

\begin{definition} A point $x\in [0,1]$ is called {\em non-wandering with respect
to  $A $} if,
 for each $ \epsilon > 0
$, there exists $k\geq 1$ and a $k$-path $ (x_1,....,x_k)$ in $
\mathcal P_k ( x,  x) $ such that
$$ \left |\sum_{i= 1}^{k - 1} (A - m) ( x_i,  x_{i+1}) \right | < \epsilon. $$
We will denote by $ \Omega(A) $ the set of non-wandering points with
respect to $ A$.

\end{definition}

The above definition is analogous (to the case of discrete time dynamics) to the continuous time one in Mather Theory (see \cite{Fathi}, \cite {CI} and \cite{GLM}). There, a point  $x$ is non wandering for the Lagrangian $L$, if you can move from $x$ to $x$ by means of connecting paths $\gamma$, with action $\int L(\gamma, \gamma' ) \, dt$  so small as you want.

\begin{lemma}\label{le:projsup} Suppose that the observable $A$ is $C^2$, and $\frac{\partial^2 A}{\partial x\partial y}>0$. Let $\nu\in\mathfrak M_{0}$ be any maximizing measure. We claim that  $\pi_1(\supp(\nu))\subset\Omega(A).$ \end{lemma}
\textit{Proof:} Let $u$ be a backward calibrated subaction, and $\mbox{dom}( Du)$ be the set of differentiable points of  $u$.
Let $Y_0:\mbox{dom}( Du)\to[0,1]$ be the map defined by  $Y_0(x)=y$, where $y$ is the unique point such that   $(x,y)$ satisfies (\ref{b}) (see the remark after theorem \ref{grafico}). As we will see in proposition
\ref{monotonicidade}, this map is monotonous, hence we can define a measurable map  $
Y: [0,1] \to [0,1]$, by
$Y(x)=Y_0(x)$ if $x \in \mbox{dom}( Du)$, and
$Y(x)=\lim_{z\rightarrow x^-,z \in \mbox{dom}( Du)}Y_0(z)$ if $x \notin \mbox{dom}( Du)$.
 Note that  $\nu_{\infty} \circ \pi_1^{-1}$-a.e. $\pi_1(\supp\nu_{\infty})\subset\mbox{dom}( Du)$.

Let us prove that $\nu_{\infty} \circ \pi_1^{-1}$ is an invariant measure
  for  $Y$.
Indeed, for $f\in C^0(\Omega(A))$, we have that:

$$\int f\circ Y(x) \,d \nu_{\infty} \circ \pi_1^{-1} (x) =
\int f\circ Y(x) \,d \nu_{\infty} (x,y)=
\int f(y) \,d \nu_{\infty} (x,y)=$$
$$=\int f(x) \,d \nu_{\infty} (x,y)=
\int f(x) \,d \nu_{\infty} \circ \pi_1^{-1} (x) \, ,
$$ where in the second equality we  used the fact that, if $(x,y)\in\supp(\nu_{\infty})$, then $y=Y(x)$,  and in the third equality we  used item (a) of lemma \ref{holon}.

Take $(x,y)\in\supp \nu_{\infty}$ and $B$ a ball centered in  $x$. We can see that $\pi_1^{-1}(B)$ is an open set which contains $(x,y)$, and this implies $\nu_{\infty} \circ \pi_1^{-1}(B)>0$. Using Poincar\'e recurrence theorem, there  exists $x_1\in B\cap\mbox{dom}( Du)$ such that, for infinitely many $j's$, $x_{j+1}:=Y^j(x_1)$ is in $B$.

Note that the points  $x_j$ satisfy the following equation:
$$u(x_j)-u(x_{j+1})=A(x_j,x_{j+1})-m $$ because, by lemma \ref{subdif},  $u$ is differentiable in each  $x_j$ and then there exists only one  $y(x_j)$ (that coincides with $x_{j+1}$) which satisfies the equation   (\ref{b}).

We fix  $\epsilon>0$ and $x_j\in B$, we can construct the following path: $(\tilde x_1,...,\tilde x_j)=(x, x_2,...,x_{j-1}, x)$, and we have that
$$\sum_{i=1}^{j-1}(A-m)(\tilde x_i,\tilde x_{i+1})=$$
$$=u(x_1)-u(x_j)+A(x,x_2)-A(x_1,x_2)+A(x_{j-1},x)- A(x_{j-1},x_j) \leq \epsilon $$ if  $B$ is small enough, because $u$ is Lipschitz  (and $A$ is $C^2$).\cqd

\begin{definition}Let us define $$ S_k(x, y) =
\inf_{ (x_1,....,x_k) \in  \mathcal P_k (x, y)  } \left [- \sum_{i =
1}^{k - 1} (A - m  )( x_i,  x_{i+1}) \right]. $$
 We call  {\em Ma\~n\'e potential}  the function $ S :  [0,1] \times [0,1] \to
\mathbb R  $ defined
by
$$ S (x, y) = \inf_{k}  S_k( x, y), $$
and {\em  Peierls barrier}  the function $h:[0,1] \times [0,1] \to
\mathbb R \cup \{ + \infty \} $ defined by
$$ h (x, y) = \liminf_{k\to\infty }  S_k( x, y). $$

\end{definition}

The value $h(x,y)$ (or, $S(x,y)$) measures, in a certain sense, the cost to move from $x$ to $y$. This will be a main tool for showing the uniqueness of the calibrated subaction.

\vspace{.4cm}

It is easy to see that
$$\Omega(A)=\{x\in[0,1] \; :\;S(x,x)= h(x,x)=0\}\,.$$

The functions $S$ and $h$ have the following  properties

(a) if $x,y,z\in[0,1]$ then $S(x,z)\leq S(x,y)+S(y,z).$

(b)  $S(\cdot, y)$ is a forward-subaction and $S(x,\cdot)$ is a backward-subaction.

(c)   $h(\cdot,y)$ is a calibrated forward-subaction and $h(x,\cdot)$ is a calibrated backward-subaction.

\begin{definition} We will say that a property  is {\em generic    for  $A$, $A\in C^{2}([0,1]^2)$, in Ma\~n\'e's sense}, if
the property is true for $A+f$, for any $f$, $f\in C^{2}([0,1])$, in a set $G$ which is generic (in Baire sense).


\end{definition}

We want to prove that, for $A$ which is generic in  Ma\~n\'e's sense \cite{Man},  the functions $V$ and $\bar V$ are unique (up to a constant). To do that, first we show that generically the maximizing measure is unique, as we will see in the following proposition.

\begin{proposition} \label{generic}Suppose that the observable $A$ is $C^2$, and $\frac{\partial^2 A}{\partial x\partial y}>0$. Then the set $$ G_2=\{f\in C^{2}([0,1])\;\;|\;\;\mathfrak M_{0}(A+f)=\{\nu\} \mbox{ and }\pi_1(supp(\nu))=\Omega(A+f)\}$$ is generic (in Baire sense) in $C^{2}([0,1])$.\end{proposition}

We will use a result  of \cite{BC} in order to prove proposition \ref{generic}.
First we will show that

\begin{equation}\label{G1} G_1=\{f\in C^{2}([0,1])\;\;|\;\;\mathfrak M_{0}(A+f)=\{\nu\}\} \end{equation} is generic (in Baire sense).
\vspace{0.2cm}

{\bf Remark :} We point out that if one considers above, in the definition of $G_2$, potentials of the form $A(x,y) \,+ \,l\, \,x$, where $l$ is constant,
instead of $A(x,y) + f(x)$, the same result is true for a generic $l\in \mathbb{R}$. This new statement is natural (and means something interesting) once it is common to consider a magnetization as a function of this form.  In this way, for example, considering fixed the term  $\frac{1}{2} (x-y)^2$, for a dense set of $l\in \mathbb{R}$, we have that the zero-temperature state for $ A(x,y)=\frac{1}{2} (x-y)^2+ l\, x$ is unique.
\vspace{0.1cm}

Let us fix some notation : $C$ is the set of continuous functions in $[0,1]^2$, $F=C^*$  the vector space of continuous functionals $\nu:C\to \re$, $E=C^{2}([0,1])$ provided with the $ C^{2}$ topology, and $G$ is the vector space of finite Borel signed measures on $[0,1]$.
$K\subset G$ is the set of Borel probability measure on $[0,1]$, and note that $\mathbf{M}\subset F$. We denote by $F_A:\mathbf{M}\to \re$  the linear functional defined by
$F_A(\nu)=\int A\;d \nu$. Note that $\mathfrak M_{0}(A)$  is the set of points of $\nu\in\mathbf{M}$ which maximize $F_A|_{\mathbf{M}}$. Finally, let  $\pi:F\to G$ be the projection  induced by  $\pi_1:[0,1]^2\to [0,1]$.

\begin{lemma} \label{BC}There exists a generic subset $\mathcal O\subset E$ (in Baire sense) such that, for all $f\in\mathcal O$, we have $$\#\pi( \mathfrak M_{0}(A+f))= 1.$$
\end{lemma}
\textit{Proof:} We just note that  $F_A$ is a affine subspace of dimension 0 of ${\mathbf{M}}^*$, then proposition follows by theorem 5 of \cite{BC}.\cqd

Note that, in order to have (\ref{G1}), we need to prove that $\# \mathfrak M_{0}(A+f)= 1$.
\begin{lemma} \label{dimension}If the observable $A$ is $C^2$, and $\frac{\partial^2 A}{\partial x\partial y}>0$, then we have $\# \mathfrak M_{0}(A)=\# \pi( \mathfrak M_{0}(A))$
\end{lemma}
\textit{Proof:} 
By lemma \ref{uniaosuportes} we know that the restriction to
$\cup_{\nu \in \mathfrak M_{0}}  \supp(\nu)$
of the projection $[0,1]^2\to [0,1]$ is a injective map. Hence the linear map $\pi: \mathfrak M_{0}(A)\to G$ is injective, and $\# \pi( \mathfrak M_{0}(A))=\# \mathfrak M_{0}(A)$.\cqd

\vspace{.5cm}
\textit{Proof of proposition \ref{generic}:} Note that, by lemmas \ref{BC} and \ref{dimension}, we have that the set $G_1$ given in (\ref{G1}) is  generic.

Let $f_0\in  G_{1}, $ and $f_1\in C^{2}([0,1])$ such that $f_1\geq 0$ and $\{x: f_1(x)=0\}=\pi_1(\supp(\nu))$. Then $\pi_1(\supp(\nu))\subset \Omega(A+f_0+f_1)$.

Claim: If $x_1\notin\pi_1(\supp(\nu))$ then $x_1\notin\Omega(A+f_0+f_1)$.

Indeed, $f_1(x_1)>0$, and
 $$h^{(A+f_0+f_1)}(x_1,x_1)=\liminf _{k\to\infty}\left(\inf_{\mathcal P_k(x_1,x_1)}\sum _{i=1}^{k-1} (A+f_0+f_1-m)(x_i,x_{i+1})    \right)\geq$$

$$ \liminf _{k\to\infty}\left(\inf_{\mathcal P_k(x_1,x_1)}\sum _{i=1}^{k-1}(A+f_0-m)(x_i,x_{i+1})  +f_1(x_1)  \right)=
$$
$$h^{(A+f_0)}(x_1,x_1)+f_1(x_1)>0.$$
Hence $\pi_1(\supp(\nu))= \Omega(A+f_0+f_1).$\cqd
\vspace{.5cm}


\begin{proposition}\label{sup}
If $u$ is a calibrated backward-subaction, then for any $x$ we have
$$
u(x) = \sup_{p \in \Omega(A) }\{u(p) - h(p,x) \}.
$$
\end{proposition}

\textit{Proof}: 
For $(x_1,...,x_k)\in\mathcal P_k( x,\bar x)$, we have $$u(x_i)-u(x_{i+1})\geq A(x_i,x_{i+1})-m$$
 and $$u(x_k)-u(x_{1})\leq -\sum_{i=1}^{k-1}A(x_i,x_{i+1})-m\,.$$ Hence,
 $ u( \bar x)-u( x)\leq h( x, \bar x)$, and therefore  $$\displaystyle
u(x) \geq\sup_{p \in \Omega(A) }\{u(p) - h(x,p) \}.$$

 Now we show the other inequality. We denote by $x_1=x$. The fact that $u$ is a backward calibrated subaction implies the existence of $x_2$ such that  $u(x_1)=u(x_{2})+ A(x_1,x_{2})-m$. Thus, recursively, we can construct $(x_1,x_2,...,x_n,...)$ such that  $u(x_n)=u(x_{n+1})+ A(x_n,x_{n+1})-m$.

  Let  $p$ be an accumulation point of the sequence  $\{x_n\}$. We claim that  $p\in\Omega(A)$. Indeed, if $x_{n_j}\to p$, we fix $j>i$, and then we construct  $(\tilde x_1,...,\tilde x_{n_j-n_i})=(p,x_{n_{i+1}},...,x_{n_{j-i}},p)$.  Hence, we have

  $$\sum _{i=1}^{n_j-n_i-1}(A-m)(\tilde x_i,\tilde x_{i+1})=$$
  $$ =\sum _{k=n_i}^{n_{j-1}}(A-m)( x_k, x_{k+1})+A(p,x_{n_{i+1}})-A(x_{n_{i}},x_{n_{i+1}})+A(x_{n_{j-1}},p)-A(x_{n_{j-1}},x_{n_{j}})=$$

 $$=u(x_{n_j})-u(x_{n_i})+A(p,x_{n_{i+1}})-A(x_{n_{i}},x_{n_{i+1}})+A(x_{n_{j-1}},p)-A(x_{n_{j-1}},x_{n_{j}})$$

 Then for  $\epsilon>0$ fixed and $i$ large enough we have that
 $$\left|\sum _{i=1}^{n_j-n_i-1}(A-m)(\tilde x_i,\tilde x_{i+1})\right|\leq \epsilon\,.$$ Therefore $p\in\Omega(A)$.

 Now take  $(\tilde x_1,...,\tilde x_{n_j})=( x_1,x_2..., x_{n_j-1},p)$. We have
 $$ -\sum _{i=1}^{n_j-1}(A-m)(\tilde x_i,\tilde x_{i+1})+u(x)-u(p)=$$
 $$= -\sum _{i=1}^{n_j-1}(A-m)( x_i, x_{i+1})+A(x_{n_j-1},x_{n_j})-A(x_{n_j-1},p)+u(x)-u(p)= $$
 $$=u( x_{n_j})-u(p)+A(x_{n_j-1},x_{n_j})-A(x_{n_j-1},p)\,.$$ Given  $k>0$ there exists  $n_k$ such that
 $$ -\sum _{i=1}^{n_k-1}(A-m)(\tilde x_i,\tilde x_{i+1})\leq u(p)-u(x)+\frac{1}{k}\,.$$ Making  $k\to\infty$ we obtain  $h(x,p)\leq u(p)-u(x) $. Then
 $$u(x)= \sup_{p\in \Omega(A)}\{ u(p)-h(x,p)\}. $$\cqd

\begin{proposition}\label{bijecao}There exists a bijective correspondence between
the set of calibrated backward-subactions and the set of functions $f\in
C^0(\Omega(A))$ satisfying $f( y)-f(x)\leq h(x,y)$,
 for all points $x,y$ in $\Omega(A)$.
\end{proposition}
\textit{Proof:} Let us suppose that $f$ satisfies $f(y)-f(x)\leq h(x,y)$. We define the following map
$\displaystyle f\mapsto u_f(x):=\sup_{p\in \Omega(A)}\{ f(p)-h(x,p)\} $. We will just show that this map is a bijection. The proof of the fact that $u_f$ is a calibrated backward-subaction is similar to the proof of theorem 13 in \cite{GL}.

We will prove that the map is injective: let  $f\in C^0(\Omega(A))$  satisfying   $f(y)-f(x)\leq h(x,y)$. For $x\in\Omega(A)$, we have that $h(x,x)=0$, and hence
$$f(p)-h(x,p)\leq f(x)\leq \sup_{p\in \Omega(A)}\{ f(p)-h(x,p)\}=u_f(x)\,.$$ Then $u_f(x)=f(x), \forall x\in\Omega(A)$. Therefore $f\neq \tilde f$ implies $u_f\neq u_{\tilde f}$.

Now, we will prove that the map is surjective: let  $u$ be a calibrated subaction. Define  $f=u|_{\Omega(A)}$. By proposition  \ref{sup}, we have that  $f$ satisfies
$f(y)-f(x)\leq h(x,y)\;\;$ and $\;\;\displaystyle u(x)=\sup_{p\in \Omega(A)}\{ u(p)-h(x,p)\}=\sup_{p\in \Omega(A)}\{ f(p)-h(x,p)\}=u_f(x).$

 \cqd

\bigskip
Now suppose that $A$ has a unique maximizing measure $\nu_{\infty}$ and also that $\pi_1 (\supp(\nu_{\infty})) = \Omega(A)$.
As we have explained in the proof of  lemma \ref{le:projsup} above, we can define a measurable map $ Y: \Omega(A) \to \Omega(A)$. Indeed, when $x$ is such that  there is unique $y$ satisfying $(x,y)\in \supp(\nu_{\infty})$, then $y=Y(x)$. In the other case, we define $Y$ via the limit coming  from the left side.

\begin{lemma} If $A$ is generic in the Ma\~n\'e sense, then
the measure $\nu_{\infty} \circ \pi_1^{-1}$ is an invariant ergodic measure for $Y$.
\end{lemma}
\textit{Proof:} First we prove the invariance:
Let $f\in C^0(\Omega(A))$.
We have:

$$\int f\circ Y(x) \,d \nu_{\infty} \circ \pi_1^{-1} (x) =
\int f\circ Y(x) \,d \nu_{\infty} (x,y)=
\int f(y) \,d \nu_{\infty} (x,y)=$$
$$=\int f(x) \,d \nu_{\infty} (x,y)=
\int f(x) \,d \nu_{\infty} \circ \pi_1^{-1} (x) \, .
$$

Now we will prove that $Y$ is uniquely ergodic: let $\eta$ be a measure in the Borel sets of $\Omega(A)$ which is invariant for $Y$.
If we define, for each Borel set $B$ of $[0,1]^2$,  $\nu(B)= \eta(\pi_1 (B\cap \supp(\nu_{\infty})))$, we have that $\nu$ is a measure probability in
$[0,1]^2$ such that

(1) $\supp(\nu) \subset \supp(\nu_{\infty})$,

(2) $\pi_1(\nu)=\eta$,

(3) $\nu \in  \mathbf{M}$.

In order to prove (3), consider $f\in C([0,1])$. We have
$$\int f(y) d\nu(x,y)= \int f(Y(x))d\nu(x,y)= \int f(Y(x))d\nu(x) =
$$
$$= \int f(x) d\nu(x)= \int f(x)d\nu(x,y)\, ,
$$
where we used, in sequence: (1); (2) ; $\eta$ is $Y$-invariant ; (2).

Note that for any calibrated backward-subaction $u$ we have
$$ \int A(x,y) d\nu(x,y) = \int \left( u(x)-u(y)+m  \right) d_\nu(x,y) = m \, ,$$ where in the first equality we used
(1) and proposition \ref{cohomologia}, and, in the second equality we used (3). Thus we have that $\nu$ is a maximizing measure, and by uniqueness $\nu=\nu_{\infty}$. This implies $\eta=\pi_1(\nu_{\infty})$, which shows that there exists an unique invariant measure for $Y$, which is a ergodic measure. \cqd

\begin{proposition}\label{ergodic} If $\nu\circ\pi_1^{-1}$ is an ergodic measure in
$[0,1]$, and $u,u'$ are two calibrated backward-subactions for $A$, then
$u-u'$ is constant in $\pi_1(\supp(\nu))$.
\end{proposition}
For the proof of this proposition see theorem 17 of \cite{GL}.

\begin{theorem} \label{unitary} If $A$ is generic in the Ma\~n\'e sense, then the set of calibrated backward-subactions has an unique element.
\end{theorem}
\textit{Proof:} By the hypothesis  $\nu_{\infty}$ is the unique
maximizing  measure, hence  $\nu_{\infty}\circ\pi_1^{-1}$ is
ergodic, and
$\pi_1(\supp(\nu_{\infty}))=\Omega(A).$

Let $f,f':\Omega(A)\to \re$ be continuous functions satisfying the
hypothesis of proposition \ref{bijecao}. In the proof of  proposition \ref{bijecao} we see that
we can get two calibrated subactions $u_f,\,u_{f'}$ such that
$f-f'=u_f-u_{f'}$ in $\Omega(A)$, and hence, by proposition \ref{ergodic}
$u_f-u_{f'}$ is constant in $\Omega(A)$. Again, from proposition
\ref{bijecao}, we show that the set of calibrated backward-subactions has an unique element. \cqd

If we consider $V$ and $\bar V$ given in lemma \ref{subaction2},
theorem \ref{unitary} proves that $\bar V$ is unique. The proof that $V$ is unique uses similar
arguments .


\section{The shift in the  Bernoulli space $[0,1]^\mathbb{N}$, and a Large Deviation Principle}\label{sec_ldp}

Let us come back to the maximization problem, over  $\mathcal M_0$, of
\begin{equation}\label{maxim2} \int A
  d\mu \,. \end{equation}

We get in this section (and from what we proved before)  a family of absolutely continuous Markov measures  $\mu_{\beta}$, indexed by a real parameter $\beta$, and this family  of measures weakly  converges, when $\beta \to \infty$, to the maximizing measure $\mu_{\infty}$.  A natural question is to know the speed (in logarithm scale)  of  convergence of the probability $\mu_{\beta}(C)\to 0$, of a $\mu_{\infty}$-null set $C$, when $\beta \to \infty$. In this direction we will present a Large Deviation Principle. This is our main goal in this section.

The following proposition allows us to conclude that, generically in Ma\~n\'e's sense, all such
maximizing measures, after projection in the first two coordinates,
are unique.

\begin{proposition}\label{medidamax}
Suppose that $\nu_{\infty}$
 is a   maximizing measure in  $\mathbf{M}$ .

(i) If $A$ has an unique maximizing measure in  $\mathbf{M}$, then any maximizing measure in   $\mathcal M_0$ is projected by
$\Pi$ in  $\nu_{\infty}$, where $\Pi:\espaco \to [0,1]^2$ is the projection in the first two coordinates.

\medskip

(ii) $\nu_{\infty}$ can be extended to a maximizing measure
$\mu_{\infty}\in \mathcal M_0$  which is a stationary Markov  measure.

\medskip

(iii) If  $\nu_{\beta}$ is the  family of measures given by
(\ref{nubeta}), then this measures can be extended to absolutely continuous Markov measures  $\mu_{\beta}$, and this sequence of measures weakly  converge to the maximizing measure $\mu_{\infty}$.
\end{proposition}

\textit{Proof}: Item  (i) follows by items  (b) and   (c) of  proposition   \ref{relacao} and
by proposition  \ref{generic}.  Item (ii) follows by  item  (a) of  proposition
\ref{relacao}. Item  (iii) is a consequence of the remark after the proof of      proposition  \ref{relacao}. \cqd

\medskip

 From now on, until the end of this section, we will suppose that the maximizing measure $\nu_{\infty}$, and the functions $V$ and $\bar V$ are unique. This is a generic property in Ma\~n\'e sense.

Thus, for the maximization problem in the Bernoulli shift,  we have shown the existence of a maximizing measure $\mu_{\infty}$ which can be approximated by absolutely continuous  stationary Markov measures $\mu_{\beta}$, which were explicitly calculated.

Now we will show a Large Deviation Principle for the  family of measures
$\{\mu_{\beta}\}$. We will also exhibit a Large Deviation Principle for the bidimensional measures $\nu_{\beta}$ which, by the earlier sections, converge to $\nu_{\infty}$.

\vspace{.6cm}

\medskip

\begin{lemma}\label{bbb} Suppose $k\geq 2$. Let $F_{k}:[0,1]^{k}\to\re$ be the function given  by
$$F_k(x_1,..,x_{k}):=\max(V+\bar V)-V(x_1)-\bar V(x_k)-\sum_{i=1}^{k-1}(A-m)(x_i,x_{i+1}).$$
Let $D_{k}=A_1....A_k$ be a cylinder of size $k$.
Then, there exists the limit
$$\lim_{\beta\to\infty}\frac{1}{\beta}\log
\mu_{\beta}(D_{k}) =  -\inf_{(x_1,..,x_{k})\in D_{k}}F_k(x_1,..,x_{k}).$$

\end{lemma}

\textit{Proof: }
Let us define
$$f_{k,\beta}(x_1,...,x_k):=$$
$$:=\frac{1}{\beta}\log \pi_\beta+\frac{k-1}{\beta}\log \lambda_\beta-\sum_{i=1}^{k-1} A(x_i,x_{i+1})-\frac{1}{\beta}\log \varphi_\beta(x_1)-\frac{1}{\beta}\log \bar\varphi_\beta(x_k).
$$
We have that $f_{k,\beta} \to F_k$ uniformly when $\beta \to \infty$.  This is a consequence of the uniqueness of $V$ and $\tilde V$.

We begin by proving the

{\it Claim}: Let $C_{k}=A_1....A_k$ be a cylinder  of size $k$. We have
$$\limsup_{\beta\to\infty}\frac{1}{\beta}\log
\mu_{\beta}(C_{k})\leq -\inf_{(x_1,..,x_{k})\in C_{k}}F_k(x_1,..,x_{k})\,.$$

To prove the {\it Claim}, note that we have
$$\mu_{\beta}(C_{k})=$$
$$=     \int_{A_1}...\int_{A_k}\frac{e^{\beta A(x_{k-1},x_k)}\bar \varphi_{\beta}(x_k)}{\bar \varphi_{\beta}(x_{k-1})\lambda_{\beta}}\;...\frac{e^{\beta A(x_1,x_2)}\bar \varphi_{\beta}(x_2)}{\bar \varphi_{\beta}(x_1)\lambda_{\beta}}\;\frac{\varphi_{\beta(x_1)}\bar\varphi_{\beta(x_1)}}{\pi_{\beta}}dx_k...dx_1=$$

\begin{equation}\label{medcil}=\int_{A_1}...\int_{A_k} e^{-\beta f_{k,\beta}(x_1,...x_k)}dx_k...dx_1\leq e^{-\beta\inf_{C_k}f_{k,\beta}(x_1,...x_k)}|C_k|\,,
\end{equation}
where $|C_k|$ denotes the Lebesgue measure of $C_k$.
Hence
$$ \frac{1}{\beta}\log \mu_{\beta}(C_{k})\leq -\inf_{C_k}f_{k,\beta}(x_1,...x_k)+\frac{1}{\beta}\log |C_k|, $$
and then, by the uniform convergence, we have:

$$ \limsup_{\beta\to\infty}\frac{1}{\beta}\log \mu_{\beta}(C_{k})\leq -\inf_{C_k}F_k(x_1,...x_k)\,,$$
which finishes the proof of the {\it Claim}.

Now we will prove the lemma:
if we fix $\delta>0$, using the continuity of $F_k$ we can find a point $(x_1,...,x_k) \in D_k^0$ (the interior of $D_k$) such that
\begin{equation}\label{aberto}\inf_{D_k}F_k\leq F_k(x_1,...,x_k)< \inf_{D_k}F_k+\delta\;\; .\;\;  
\end{equation}
Now, let $D_{\delta}$ be a  cylinder of size $k$, such that $(x_1,...,x_k) \in D_{\delta} \subset D_k^0$, and
\begin{equation}\label{aberto2}
\inf_{D_k}F_k\leq F_k(y_1,...,y_k)< \inf_{D_k}F_k+2 \delta\;\;\;  \forall (y_1,...,y_k)\in D_{\delta}\,.
\end{equation}
We have that

$$\mu_{\beta}(D_{k})\geq \mu_{\beta}(D_{\delta})\geq e^{- \beta \sup_{D_{\delta}} f_{k,\beta}(y_1,...,y_k)} |D_{\delta}|\;\;,$$
where the last inequality cames from (\ref{medcil}).
Now we use again the uniform convergence of  $f_{k,\beta}$ to $ F_k $ in order to get

$$\liminf_{\beta\to \infty}\frac{1}{\beta}\log\mu_{\beta}(D_{k})\geq -\sup_{D_{\delta}}F_k\,.$$

By (\ref{aberto2}), we get

\begin{equation}\label{jjj}
\liminf_{\beta\to\infty}\frac{1}{\beta}\log\mu_{\beta}(D_{k})\geq-\inf_{D_k}F_k- 2\delta
\end{equation}
Sending $\delta \to 0$, and using the {\it Claim}, we finish the proof of the lemma.
\cqd
\vspace{.2cm}

Note that if we set $k=2$ above, we get a LDP for the family $\nu_{\beta} \to \nu_{\infty}$.

\begin{theorem}\label{ldp} Let $I:[0,1]^{\nat}\to [0,+\infty]$ be the function defined by
$$I(\mathbf{x}):= \sum_{i\geq 1}V(x_{i+1})-V(x_i)-(A-m)(x_i,x_{i+1})\,.   $$

 Let $D=A_1....A_k$ be a cylinder of any size $k$.
Then, there exists the limit

$$\lim_{\beta\to\infty}\frac{1}{\beta}\log
\mu_{\beta}(D)= -\inf_{\mathbf{x}\in D}I(\mathbf{x})\,.$$


\end{theorem}


\vspace{.2cm}

Note that, by lemma \ref{subaction2}, $V(x_{i+1})-V(x_i)-A(x_i,x_{i+1})+m\geq 0$, therefore  the sequence of partial sums of the series in the definition  of $I(\mathbf{x})$ is a non-decreasing sequence. This shows that $I(\mathbf{x})$ is well defined (note that $I(x)$ can be $+\infty$).

In order to prove Theorem \ref{ldp} we will need some new results and definitions.

For each $N\geq 2$, let us extend the function $F_N$ to the space
$\espaco$:
$$F_N (\mathbf{z}) :=F_N(z_1,...,z_n)= \sup(V+\bar V) - V(z_1)-\bar V(z_{N}) - \sum_{i=1}^{N-1} (A-m)(z_i,z_{i+1}) \;.$$

\begin{lemma}\label{Desig_F_N}
$\forall \;\mathbf{z} \in [0,1]^{\mathbb{N}}$, 
 we have
$$F_N(\mathbf{z}) \geq \max(V+\bar V) - ( V(z_1)+\bar V(z_{1})) \geq 0 \;. $$
\end{lemma}

\textit{Proof: }
By lemma \ref{subaction2}

$$\bar V(x)-\bar V(y)\geq A(x,y)-m\;\;,\; \forall x,y,$$
then

$$-\sum_{i=1}^{N-1}(A-m)(z_i,z_{i+1})\geq \bar V(z_{N})-\bar V(z_1)\,.$$
Hence, by definition of $F_N$
$$F_N (\mathbf{z})  \geq \max(V+\bar V) - ( V(z_1)+\bar V(z_{1}))\,.$$

\cqd
\begin{lemma}\label{Lim_V_Vestrela_Existe}
(a)  for a fixed $\mathbf{x}\in \espaco$, we have that
$$V(x_1)+\sum_{i=1}^{k-1}(A-m)(x_i,x_{i+1}) +\bar V(x_k)$$ is decreasing with respect to $k$.

\medskip
(b) If $I(\mathbf{x})<+\infty$, then there exists the limit
$$L(\mathbf{x})= \lim_{k \to +\infty} V(\sigma^k(\mathbf{x})))+\bar V(\sigma^k(\mathbf{x}))).$$
\end{lemma}

\textit{Proof: }
(a)  $$V(x_1)+\sum_{i=1}^{k}(A-m)(x_i,x_{i+1}) +\bar V(x_{k+1})=$$ $$=V(x_1)+\sum_{i=1}^{k-1}(A-m)(x_i,x_{i+1}) +\bar V(x_{k})+A(x_k,x_{k+1})-m+\bar V(x_{k+1})-\bar V(x_{k}),$$
by lemma \ref{subaction2} (remember that $\tilde m=m$) $\;A(x_k,x_{k+1})-m+\bar V(x_{k+1})-\bar V(x_{k})\leq 0$, and we have (a).

\medskip
(b) We have
$$I(\mathbf{x})=\sum_{i\geq 1}V(x_{i+1})-V(x_i)-(A-m)(x_i,x_{i+1}) =$$
\begin{equation}\label{I}=\lim_{k \to +\infty} V(x_k)+\bar V(x_k)-\lim_{k\to\infty}
\left(V(x_1)+\sum_{i=1}^{k-1}(A-m)(x_i,x_{i+1}) +\bar V(x_k)\right) \end{equation}

\medskip

Hence, if $I(\mathbf{x})<+\infty$, it follows,  thanks to item (a), that  $ V(x_k)+\bar V(x_k)=V(\sigma^k(\mathbf{x})))+\bar V(\sigma^k(\mathbf{x})) $ must converge.

\cqd


\begin{lemma}\label{aaa}
Suppose $I(\mathbf{x})<+\infty$. 
Then, if we define, for each $M \in \mathbb{N}$, the probability measure
$$ \mu_M = \frac{1}{M} \sum_{j=1}^{M-1}\delta_{\sigma^j(\mathbf{x})} \;, $$
we have that $\Pi(\mu_M) \to \nu_{\infty}$ in the weak-$\star$ topology (where $\Pi$ is the projection in the two first coordinates).

\end{lemma}

\noindent {\textit{Proof: }
Given $\epsilon >
0$, there exists $N_{\epsilon} \in \mathbb{N}$ such that , for all
$N \geq N_{\epsilon}$, and all $M > N$,
$$\sum_{i = N}^{M-1} V(x_{i+1})-V(x_i)-(A-m)(x_i,x_{i+1}) < \epsilon \;.$$
Thus
$$
V(x_M)-V(x_N)+(M-N) m  <  \sum_{i = N}^{M-1} A(\sigma^i (\mathbf{x}) ) + \epsilon \;,
$$
and
$$
\frac{1}{M-N}\sum_{i = N}^{M-1} A(\sigma^i (\mathbf{x}) )   >   m + \frac{V(x_M)-V(x_N)}{M-N}  - \frac{\epsilon}{M-N} \;,$$
and then we get that
$$ 
\liminf_{M\to +\infty}  \frac{1}{M}\sum_{i = 1}^{M-1} A(\sigma^i (\mathbf{x}) )   \geq   m \;.$$
Now we remember that
$$ \frac{1}{M}\sum_{i = 1}^{M-1} A(\sigma^i (\mathbf{x}) ) = \int A d\mu_M  \leq m\;, $$
and finally we  get
$$ \lim_{M\to +\infty} \int A d\mu_M =
\lim_{M\to +\infty}  \frac{1}{M}\sum_{i = 1}^{M-1} A(\sigma^i (\mathbf{x}) )   =  m \;.$$

If we use the compactness of the  closed ball of radius 1 in the   weak-$\star$ topology, we get that $\{ \mu_M \}$ has convergent subsequences. Any limit of a convergent subsequence is a stationary measure (a $\sigma$-invariant measure) and must be a maximizing measure,
by the last equality.
As any maximizing measure is projected in $\mu_{\infty}$ by $\Pi$, we get  the lemma.
\cqd

\vspace{.3cm}


\begin{proposition}\label{Lim_V_Vestrela}
If $I(\mathbf{x})<+\infty$, then
$$\lim_{k \to +\infty} V(\sigma^k(\mathbf{x}))+\bar V(\sigma^k(\mathbf{x})) = \max (V+\bar V). $$
\end{proposition}

\noindent {\textit{Proof: }
Let ${\mathbf{z}} =(z_1,z_2,z_3,...) \in \supp(\mu_{\infty})$. We have that $(z_1,z_2) \in \supp(\nu_{\infty})$. Thus, by lemma \ref{aaa}  there exists a sub-sequence such that $\Pi(\sigma^{k_l}(\mathbf{x})) \to (z_1,z_2)$.

Fix $\epsilon > 0$.
Let
$$B_{k_l,\epsilon}(\mathbf{x}) := \{ \mathbf{y} \in \espaco \; : \; |y_j - x_{j+k_l} | \leq \epsilon, \; \forall \; 1 \leq j \leq 2 \}\;$$
be the closed cylinder of size $2$ 'centered' at
$\sigma^{k_l}(\mathbf{x})$.

If $l$ is big enough, we have that
$$
B_{k_l,\epsilon}(\mathbf{x}) \subset \{ \mathbf{y} \in \espaco \; : \; |y_j -  z_{j} | \leq 2 \epsilon, \; \forall \; 1 \leq j \leq 2 \}\;. $$

Note that 
$\mu_{\infty}(B_{k_l,\epsilon}(\mathbf{x}))  =\nu_{\infty}(B_{k_l,\epsilon}(\mathbf{x}))  >0$,
and thus using Lemma \ref{bbb} with $k=2$, it follows  that there exists a point 
$(z_{1,\epsilon},z_{2,\epsilon},z_{3,\epsilon},z_{4,\epsilon},...) \in B_{k_l,\epsilon}(\mathbf{x})$, such that $F_2((z_{1,\epsilon},z_{2,\epsilon}))=0$.   

Then, we can  use the fact that $F_2$ depends only on its first $2$ coordinates
in order to obtain that $F_2(\mathbf{w}_{\epsilon})=0$, where
$\mathbf{w}_{\epsilon}= (z_{1,\epsilon},z_{2,\epsilon},z_3,z_4,...)$ is defined by the point of $\espaco$ whose first $2$ coordinates are equal to those of $(z_{1,\epsilon},z_{2,\epsilon})$, while the other coordinates are equal to those of $ {\mathbf{z}}$.

Now, if we send $\epsilon \to 0$, we have that $\mathbf{w}_{\epsilon} \to  {\mathbf{z}}$. Thus we can use the continuity of $F_N$ to get that $F_2( {\mathbf{z}})=0$.

Using again the continuity of $F_2$, we have that $F_2(\sigma^{k_l}(\mathbf{x})) \to 0$.

Lemma \ref{Desig_F_N}  shows that

$$\lim_{l \to +\infty} V(\sigma^{k_l}(\mathbf{x}))+\bar V(\sigma^{k_l}(\mathbf{x})) = \max (V+\bar V), $$

and finally using Lemma \ref{Lim_V_Vestrela_Existe}(b) we prove proposition \ref{Lim_V_Vestrela}.
\cqd

\vspace{.7cm}

\noindent {\textit{Proof of theorem \ref{ldp}:}}  First we need to
prove the following claim.

{\it Claim}:
$$I(\mathbf{x})=\max(V+\bar V)-\lim_{k\to\infty}\left(V(x_1)+\sum_{i=1}^{k-1}(A-m)(x_i,x_{i+1}) +\bar V(x_k)\right).  $$

In order to prove the {\it Claim}, we have to consider two possibilities:
if  $I(\mathbf{x})< +\infty$, then
(\ref{I}) can be combined with proposition \ref{Lim_V_Vestrela} to give the {\it Claim}.
If $I(\mathbf{x})= +\infty$, we just have to use the expression
$$I(x) = \lim_{k\to\infty} \left( V(x_k)-V(x_1) - \sum_{i=1}^{k-1}(A-m)(x_i,x_{i+1}) \right).$$

 Thanks to Lemma \ref{bbb},
we just have to show that
$$-\inf_{(x_1,..,x_k)\in D}F_k(x_1,..,x_k) =
-\inf_{\mathbf{x}\in D} I(\mathbf{x}) \;.$$

We begin by proving that
$$-\inf_{(x_1,..,x_k)\in D}F_k(x_1,..,x_k) \leq
-\inf_{\mathbf{x}\in D} I(\mathbf{x}) \;.$$

Given $\delta>0$,
there exists a point $(y_1,...,y_k)\in D$ such
that
$$F_k(y_1,...,y_k)<\inf_{(x_1,..,x_k)\in C}F_k(x_1,..,x_k)+\delta\,.$$
By the definition of $F_k$,
$$F_k(y_1,...,y_k)=\max(V+\bar V)-V(y_1)-\bar V(y_k)-\sum_{i=1}^{k-1}(A-m)(y_i,y_{i+1}).$$
For
each $j\geq k$ we choose a  $y_{j+1}$ that satisfies $\bar
V(y_{j})=\bar V(y_{j+1})+A(y_j,y_{j+1})- m\,.$
Then
we define $\mathbf{y}:=(y_1,...y_k,y_{k+1},...)\,.$

 \vspace{.5cm}
\textit{Second Claim:} $I(\mathbf{y})= F_k(y_1,...y_k)$. Indeed,
$$F_k(y_1,...y_k)=\max(V+\bar V)-
\left(V(y_1)+\bar
V(y_k)+\sum_{i=1}^{k-1}(A-m)(y_i,y_{i+1})\right)=$$
$$=\max(V+\bar
V)-\left(V(y_1)+\bar V(y_j)+\sum_{i=1}^{j-1}(A-m)(y_i,y_{i+1})\right),\;\;\;\;
\forall j\geq k.$$

Then, from the reasoning above and the way we choose  $\mathbf{y}$, we get that $F_k(y_1,...y_k)$ is equal to
$$\max(V+\bar
V)-\lim_{j\to\infty}\left(V(y_1)+\bar V(y_j)+\sum_{i=1}^{j-1}(A-m)(y_i,y_{i+1})\right)=I(\mathbf{y})\,.$$

This implies that $$-\inf_{(x_1,..,x_k)\in D}F_k(x_1,..,x_k) < -I(\mathbf{y}) + \delta \leq
-\inf_{\mathbf{x}\in D} I(\mathbf{x}) +\delta \,.$$
Making $\delta \to 0$, we have the first inequality.

\medskip
Now, 
we will prove the second inequality:

$$
-\inf_{\mathbf{x}\in D} I(\mathbf{x})\leq -\inf_{(x_1,..,x_k)\in D}F_k(x_1,..,x_k) \,.$$

We can use Lemma \ref{Lim_V_Vestrela_Existe}(a),
and then we get, by the \textit{Claim},
$$I(\mathbf{x})=\max(V+\bar V)-\lim_{j\to\infty}\left(V(x_1)+\sum_{i=1}^{j-1}(A-m)(x_i,x_{i+1}) +\bar V(x_j)\right)\geq$$
$$\geq \max(V+\bar V)- \left(V(x_1)+\sum_{i=1}^{k-1}(A-m)(x_i,x_{i+1}) +\bar V(x_k) \right)=F_k(x_1,..,x_k).$$
\cqd


Here, finally, we can give the proofs of  theorems  \ref{principal1} and \ref{principal}: 

\textit{Proof of theorem \ref{principal1}:} (a) It follows by proposition \ref{generic} and   item (i) of  proposition \ref{medidamax}.

(b) Theorem  \ref{unitary} shows  that, generically, the set of backward calibrated subactions has an unique element. The proof that the set of forward calibrated subactions has an unique element is similar.

\textit{Proof of theorem \ref{principal}:} (a) It follows by   items (ii) and (iii) of proposition \ref{medidamax} and theorem \ref{max}.

(b) This is theorem \ref{ldp}, note that the hypothesis are fulfilled  when theorem \ref{principal1} is true.

\bigskip

We will finish this section showing the monotonicity of the graph under the twist condition.

Suppose  $A$ is $C^2$ and satisfies $$\frac{\partial^2A}{\partial x\partial y}(x,y)>0\,.$$
Then, for all $x<x', y<y'$
we have  that \begin{equation}\label{twist}A(x,y)+A(x',y')>A(x,y')+A(x',y).\end{equation}
Let $\bar V$ be the calibrated backward-subaction define above. 

As a consequence of $A$ being $C^2$, we have that $\bar V$ is Lipschitz, hence $\bar V$ is differentiable $\lambda$-a.e., where $\lambda$ is the Lebesgue  measure. Let $\mbox{dom } (D\bar V)$ be the set of points where $\bar V$ is differentiable.

Following the proof of theorem \ref{grafico}, we have that, for $x\in \mbox{dom } (D\bar V)$, there exists only one $y(x)$ such that \begin{equation}\label{equality}\bar V(x)=A(x,y(x))+\bar V(y(x))-m.\end{equation}
\begin{proposition}\label{monotonicidade}The function $Y:\mbox{dom } (D\bar V)\to [0,1]$, defined  by $Y(x)=y(x)$, $y(x)$ satisfying (\ref{equality}), is monotone nondecreasing.
\end{proposition}
\textit{Proof: }Let $x<x'$. Let us call $z=Y(x), z'=Y(x')$, and suppose that $z>z'$. We know that
$$\bar V(x)=A(x,z)+\bar V(z)-m\;\;\;,\;\;\;\bar V(x')=A(x',z')+\bar V(z')-m,$$
and
$$\bar V(x)\geq A(x,z')+\bar V(z')-m\;\;\;,\;\;\;\bar V(x')\geq A(x',z)+\bar V(z)-m.$$

Adding the first two equation and comparing with the summation of the last two, we get that
$$A(x,z)+A(x',z')\geq A(x,z')+A(x',z),$$
for $x<x', z'<z$, which is a  contradiction with  (\ref{twist}).\cqd

\vspace{0.3cm}
If we assume that $\frac{\partial^2A}{\partial x\partial y}(x,y)<0,$ then a function  $Y(x)$ as above can be defined, and it will be monotone non-increasing.

\section{Separating subactions}\label{sect_separating}

There exist subactions which are not calibrated but that are also special. One can ask about the ones which are minimal in a certain sense: the subcohomological inequality is an equality  in the smallest possible set. This subactions are called separeted subactions. 

The main goal of this section is to show the existence of a separating
subaction (see \cite{GLT} and \cite{GLM} for related results). The
idea is: given a potential $A$, we can find a subaction $u$ such that, in the
cohomological equation, the equality just holds in points $x$ that
are on $\Omega(A)$ (where it has to hold, anyway). In this way, we
have a criteria to separate points of $\Omega(A)$ from the other
ones. We can then consider a new  potential $\tilde{A}= A(x,y)+
u(x)-u(y)$ where the maximum of $\tilde{A}$ is exactly attained in
$\Omega(\tilde A)$.

\vspace{0.5cm}

\begin{definition}
 A continuous function $u: [0,1] \to \mathbb{R} $ is called a

 (a) {\em forward-subaction}
if, for any $x,y \in [0,1]$ we have

\begin{equation}\label {fsubaction} u(y)\geq  A(x,y)+ u(x)-m.\end{equation}

(b) {\em backward-subaction} if, for any $x,y\in [0,1]$ we have
\begin{equation}\label{bsubaction} u(x)\geq A(x,y)+ u(y)-m.\end{equation}

\end{definition}



\begin{definition}
We say that a  forward subaction $u$  is {\em separating} if
$$\max_x [A(x,y)+ u(x)-u(y)]=m \iff x\in\Omega(A),$$
and a backward subaction $u$  is {\em separating} if
$$\max_y [A(x,y)+ u(y)-u(x)]=m \iff y\in\Omega(A).$$
\end{definition}

We will show the existence of a separating backward-subaction.

\begin{lemma} \label{peierl} If $x\in\Omega(A)$ there exists   $\mathbf{x}=(x_1,....,x_k,...)\in\espaco$ such that  $x_1=x$ and $${h}(x_k,x_1)\leq \sum_{i=1}^{k-1}(A-m)(x_i,x_{i+1}).$$

\end{lemma}

\textit{Proof:} If  $x\in\Omega(A)$, then there exists a sequence of  paths  $\{(x^n_1,...,x^n_{j_n})\}_{n\in\nat}$ such that $x^n_1=x^n_{j_n}=x$ and $j_n\to\infty$ satisfying
\begin{equation}\label{somatend0}\sum_{j = 1}^{j_n - 1} (A - m  )( x^n_j,  x^n_{j+1})\to 0.
\end{equation}
Because  $|x^n_j|\leq 1$,
there exists a ray  $(x_1,...,x_k,...)$ which  is the limit of the paths above,  the convergence being uniform in each compact part.

Fixed $k\in\nat$. For $j_n>k$, we have that
$$S^{j_n-k}(x_k,x_1)\leq -A(x_k,x^n_{k+1})+m-\sum^{j_n-1}_{j=k+1}(A-m)(x_j^n,x_{j+1}^n),$$
and

$$ S^{j_n-k}(x_k,x_1)+\sum_{j = 1}^{j_n - 1} (A - m  )( x^n_j,  x^n_{j+1})\leq $$
$$- A(x_k,x^n_{k+1})+m+\sum_{j = 1}^{k} (A - m  )( x^n_j,  x^n_{j+1}).$$
Hence taking the $\liminf_{n\to \infty}$ and using (\ref{somatend0}) we obtain

$${h}(x_k,x_1)\leq \sum_{j = 1}^{k-1} (A - m  )( x_j,  x_{j+1}).$$ \cqd

\begin{lemma}\label{igualdade} Let  $u$ be any backward-subaction, then for all  $x\in\Omega(A)$ we have $$\max_y \left\{u(y)-u(x)+A(x,y)\right\}= m.$$

\end{lemma}
\textit{ Proof:}
 Using the fact that $u$ satisfies  equation (\ref{bsubaction}), for any $(x_1,...,x_k)\in\mathcal P_k(x,y)$,  we have that  $u(y)-u(x)\leq -\sum_{j=1}^{k-1} (A-m)(x_i,x_{i+1})$.
Hence $u(y)-u( x)\leq {h}( x, y)$.

 Let $x\in\Omega(A)$ and let $\mathbf{x}=(x_1,...,x_k,...)$ be the point in $\espaco$  which exists by  lemma \ref{peierl}.

 By  lemma \ref{peierl} we have that
$$u(x_1)-u(x_k)\leq {h}(x_k,x_1)\leq \sum_{j = 1}^{k-1} (A - m  )( x_j,  x_{j+1}),$$


and, as it is a backward-subaction,
$$u(x_k)-u(x_1)\leq  -\sum_{j = 0}^{k-1} (A - m  )( x_j, x_{j+1}).$$

In particular, for  $k=1$, $$u(x_2)-u(x_1)=-A(x_1,x_2)+m.$$

This implies

$$\max_y \left\{u(y)-u(x)+A(x,y)\right\}= m.$$ \cqd

\begin{lemma} If the observable $A$ is H\"{o}lder continuous, then the function   $S_x(\cdot):=S(x,\cdot)$ is uniformly H\"{o}lder and has the same H\"{o}lder  constant of $A$.

\end{lemma}
\textit{Proof:} Let us fix $x$, $\epsilon>0$ and $ y, z\in[0,1]$,  then there exists $(x_1,...,x_k)\in\mathcal P_k(x,y)$ such that $$|-\sum_{i=1}^{k-1} (A-m)(x_i,x_{i+1})|\leq S(x,y)+\epsilon. $$

 Consider now the following path: $(\tilde x_1,...,\tilde x_k)=(x_1,...,x_{k-1},z)\in\mathcal P_k(x,z)$, then
$$ -\sum_{i=1}^{k-1} (A-m)(\tilde x_i,\tilde x_{i+1})=-\sum_{i=1}^{k-1} (A-m)(x_i,x_{i+1})+A(x_{k-1},y)-A(x_{k-1},z)\,.$$
Therefore,
$$S(x,z)\leq -\sum_{i=1}^{k-1} (A-m)(\tilde x_i,\tilde x_{i+1})\leq S(x,y)+\epsilon +\mbox{Hol}_{\alpha}(A) |z-y|^{\alpha},\;\;\;\; \forall \epsilon,$$
i.e., $S(x,y)-S(x,z)\leq\mbox{Hol}_{\alpha}(A) |z-y|^{\alpha}$. Changing the role of $y$ and $z$ we obtain $|S(x,y)-S(x,z)|\leq\mbox{Hol}_{\alpha}(A) |z-y|^{\alpha}$, which give us the H\"{o}lder continuity of $S_x$, independently of $x$.\cqd

\begin{theorem}\label{FS1} If the observable $A$ is H\"{o}lder continuous, there exists a separating backward-subaction.\end{theorem}

\textit{Proof:} By definition,  $$S(x,y)\leq -A(x,y)+m \;\;\;\;\forall\;y\in[0,1]\,.$$

If $x\notin \Omega(A)$, then $S(x,x)>0$. Hence
$$S_x(y)-S_x(x)<-A(x,y)+m   \;\;\;\;\forall \;y\in[0,1]\,.$$

 $\Omega(A)$ is a closed set, and thus for each $x\notin \Omega(A)$ we can find
a neighborhood $V_x\subset [0,1]\backslash\Omega(A)$ of $x$  such that
$$S_x(y)-S_x(z)<-A(z,y)+m,  \;\;\;\;\forall \;y\in[0,1], \forall\;z\in V_x\,.$$

 We can extract, from the family of these neighborhoods  $\{V_x\}_{x\notin \Omega(A)}$, a countable  family
  $\{V_{x_j}\}_{j=1}^{\infty}$ which is  a covering of $[0,1]\backslash\Omega(A)$.

  We define
  $$\tilde S_{x_j}(z)=S_{x_j}(z)-S_{x_j}(0)\,.$$

   $S_{x_j}$ is uniformly H\"{o}lder, which implies that $|\tilde S_{x_j}(z)|\leq \mbox{Hol}_{\alpha}(A) z^{\alpha},\;\;\; \forall \;x_j$, therefore the series

$$u(z)=\sum_{j=1}^{\infty}\frac{\tilde S_{x_j}(z)}{2^{j}}$$
is well defined and uniformly convergent, because $[0,1]$ is compact.
Note that $u$ is a infinite convex combination of backward-subactions $\tilde S_{x_j}$, then $u$ is also a backward-subaction.

 Fix $x\in [0,1]\backslash\Omega(A)$, there exists $k\geq 1$ such that $x\in{V_{x_k}}$. Now, $\forall y\in[0,1]$ we have
$$u(y)-u(x)=\sum_{j=1}^{\infty}\frac{S_{x_j}(y)-S_{x_j}(x)}{2^{j}}=\frac{S_{x_k}(y)-S_{x_k}(x)}{2^{k}}+\sum_{j\neq k}\frac{S_{x_j}(y)-S_{x_j}(x)}{2^{j}}<$$
$$\frac{-A(x,y)+m}{2^k}+\sum_{j\neq k}\frac{-A(x,y)+m}{2^j}<-A(x,y)+m .$$
Hence,
$$\max_y \left\{u(y)-u(x)+A(x,y)\right\}< m, \mbox{\;\;\;if\;\;\;} x\notin \Omega(A),$$
and, as $u$ is a backward-subaction, we have by lemma \ref{igualdade} that

$$\max_y\left\{u(y)-u(x)+A(x,y)\right\}= m,\;\; \mbox{\;\; if \;\;}\;\; x\in\Omega(A).$$
\cqd


\begin{thebibliography}{99}

   \vspace{0.6cm}






 \bibitem [A] {A} Adams, S. {\it  Mathematical Statistical Mechanics}. Max-Plank-Institut fur Math. (2006).




 \bibitem [AL] {AL} Athreya, K and Lahiri, S.   {\it  Measure Theory and Probability Theory}. Springer Verlag. (2006).



\bibitem [Ban] {Ban} Bangert, V. {\it Mather sets for twist maps and geodesics on tori}. Dynamics Reported. {\bf 1} (1998), 1-56.

\bibitem [BBNg] {BBNg} Brevik, I,  Borven, J-M and  Ng, S.
{\it Viscous Brane Cohomology with a Brane-Bulk energy interchange
term}. General Relativity and Gravitation. Vol. {\bf 38}, N. {\bf 5}
(2006), 907-915(9).



 \bibitem [BLT] {BLT} Baraviera, A.,  Lopes, A. O and Thieullen, Ph.
{\it A Large Deviation Principle for equilibrium states of
H\"{o}lder potentials: the zero temperature case}. Stoch. and Dyn.
{\bf 6} (2006), 77-96.


 \bibitem [Ba] {Ba} P. Bhattacharya and M. Majumdar.   {\it  Random Dynamical Systems}. Cambridge Univ. Press. (2007).

\bibitem [CG] {CG}
Chou, W. and Griffiths, R. {\it Ground states of one-dimensional
systems using effetive potentials}. Physical Review B. Vol. {\bf
34}, N {\bf 9}, (1986), 6219-6234.


\bibitem [Cv] {Cv}
Cveti, M., Nojiri, S. and Odintsov, S. D. {\it Black hole
thermodynamics and negative entropy in de Sitter and anti-de Sitter
Einstein–Gauss–Bonnet gravity}. Nuclear Physics B. Vol. {\bf 628},
Issues {\bf 1-2}, (2002), 295-330,






\bibitem [BC] {BC} Bernard, P. and Contreras, G.
{\it A Generic Property of Families of Lagrangian Systems}. Annals
of Math. Vol. {\bf 167}, No.{\bf 3}, (2008), 1099-1108



\bibitem[CI] {CI}
Contreras, G. and Iturriaga, R. {\it Global minimizers of autonomous
Lagrangians}. 22$^\circ$ Co\-l\'o\-quio Brasileiro de Matem\'atica,
IMPA, (1999).



\bibitem[CLT] {CLT}
Contreras, G., Lopes, A. O. and Thieullen, Ph. {\it Lyapunov
minimizing measures for expanding maps of the circle}. Ergodic
Theory and Dynamical Systems. Vol {\bf 21}, (2001), 1379-1409.



\bibitem [CG]     {CG}
Conze, J.P. and Guivarc'h, Y. {\it Croissance des sommes ergodiques
et principe variationnel}. manuscript circa (1993).

\bibitem [CS] {CS} Cannarsa, P. and Sinestrari, C. {\it
     Semiconcave functions, Hamilton-Jacobi equations, and optimal
     control}. Progress in Nonlinear Differential Equations and their
   Applications 58. Birkhäuser Boston Inc., Boston, MA. (2004).



 \bibitem [De] {De} Deimling, K.
{\it Nonlinear Functional Analysis}. Springer Verlag. (1985)

\bibitem [Dellach] {Dellach} Dellacherie, C.  {\it Probabilities and potential}.
North-Holland. (1978).


 \bibitem [DZ] {DZ} A. Dembo and O. Zeitouni. {\it  Large Deviations Techniques and Applications}. Springer Verlag. (1998).




\bibitem [Ev] {Ev} Evans, L. C.  {\it  Weak Convergence Methods for Nonlinear Partial Differential Equations}. Published
for the Conference Board of the Mathematical Sciences, Whashington,
DC, (1990).


\bibitem[Fathi] {Fathi}
Fathi, A. {\it Th\'eor\`eme KAM faible et th\'eorie de Mather sur
les syst\`emes lagrangiens}. Comptes Rendus de l'Acad\'emie des
Sciences, S\'erie I, Math\'ematique. Vol {\bf 324} (1997),
1043-1046.

 \bibitem [FS] {FS} Fathi, A. and Siconolfi, A. {\it Existence of   $\, C^1$ critical
 subsolutions of the Hamilton-Jacobi equations}. Inv. Math. {\bf 155} (2004), 363-388.



 \bibitem [GL] {GL} Garibaldi, E. and Lopes, A. O. {\it On Aubry-Mather theory for symbolic dynamics}.
 Ergodic Theory and Dynamical Systems.  Vol {\bf 28}, Issue {\bf 3} (2008), 791-815.


 \bibitem [GLT] {GLT} Garibaldi, E., Lopes, A. O. and Thieullen, Ph. {\it On separating sub-actions}. Preprint (2006). To appear.



\bibitem [GLM] {GLM} Gomes, D., Lopes, A. O. and Mohr, J. {\it The Mather measure and a Large Deviation Principle for the Entropy Penalized Method}.
 Preprint (2007). To appear.





 \bibitem [Go] {Go} Gole, C.   {\it  Sympletic twist maps}. World Sci. Pub Co Inc. (1998).


 \bibitem [Gom] {Gom}  Gomes, D. A. {\it Viscosity Solution methods
     and discrete Aubry-Mather problem}. Discrete Contin. Dyn. Syst.
   {\bf 13 (1)} (2005), 103-116.



 \bibitem [Gom1] {Gom1} Gomes, D. A. {\it Calculus of Variations}.
 IST - Lisboa.
 (2006).



 \bibitem [GV] {GV} Gomes, D. A. and Valdinoci, E. {\it Entropy
     Penalization Methods for Hamilton-Jacobi Equations}. Adv. Math. {\bf 215}, No. {\bf 1}, (2007), 94-152.

\bibitem[Hop] {Hop}  Hopf, E. {\it An inequality for Positive Linear Integral Operators}. Journal of Mathematics and Mechanics. Vol. {\bf 12}. N. {\bf 5} (1963),
683-692.



\bibitem[Jen1] {Jen1} Jenkinson, O. {\it Ergodic optimization}. Discrete and Continuous
Dynamical Systems, Series A. 15 (2006). 197-224.

\bibitem[Ju] {Ju} Jumarie, G. {\it Relative Information}. Springer
Verlag.
(1990).


 \bibitem [Ka] {Ka} Karlin, S.   {\it  Total Positivity}. Standford Univ. Press. (1968).

\bibitem [Lu] {Lu}
Lubkin, E. {\it Negative entropy, energy, and heat capacity in
connection with surface tension: Artifact of a model or real?}.
Inter. Journal of Theoretical Physics. Vol. {\bf 26}, N. {\bf 5}
(1987), 455-481

\bibitem[Man] {Man}
Ma\~n\'e, R. {\it  Generic properties and problems of minimizing
measures of Lagrangian systems}. Nonlinearity. Vol {\bf 9} (1996),
273-310.



\bibitem[Mat] {Mat}
Mather, J. {\it Action minimizing invariant measures for positive
definite Lagrangian Systems}. Math. Z. {\bf 207 (2)} (1991),
169-207.


\bibitem[Mo] {Mo}
Morris, I. D.  {\it A sufficient condition for the subordination
principle in ergodic optimization}. Bull. Lond. Math. Soc. {\bf 39}
no. {\bf 2}. (2007). 214-220.


\bibitem[Mi] {Mi}
Mitra, I.  {\it Introduction to dynamic optimization theory},
Optimization and Chaos. Editors M. Majumdar, T. Mitra and K.
Nishimura. Springer Verlag. (2000), 31-108.

\bibitem[Ni] {Ni} Niven, R. K.
{\it Cost of s-fold Decisions in Exact Maxwell-Boltzmann,
Bose-Einstein and Fermi-Dirac Statistics}. Physica A. Volume {\bf
365}, Issue {\bf 1} (2006), 142-149.


\bibitem[Os] {Os} Ostrowski, A. {\it On positive matrices}. Math. Annalen. Vol. {\bf 150} (1963),
276-284.

\bibitem [Pe] {Pe} Pettini, M. {\it Geometry and topology in Hamiltonian dynamics and statistical mechanics}. Springer Verlag. (2007).

\bibitem[PP] {PP}
Parry, W. and Pollicott, M. {\it Zeta functions and the periodic
orbit structure of hyperbolic dynamics}. Ast\'erisque. Vol {\bf
{187-188}} (1990).



 \bibitem [Ra] {Ra} Rachev S. and Ruschendorf, L.   {\it  Mass transportation problems, Vol I and II}. Springer Verlag. (1998).

 \bibitem [Roc] {Roc} Rockafellar, R. T.   {\it Extention of Fenchel's duality theorem for convex functions}. Duke Math. J. Vol {\bf 33} (1966)
 81-89.


\bibitem[RRS] {RRS} Risau-Gusman, S., Ribeiro-Teixeira, A. C.
and Stariolo, D. A. {\it  Topology and Phase Transitions: The Case
of the Short Range Spherical Model}. Journ. of Statist. Physics. Vol
{\bf 124} no. {\bf 5} (2006), 1231-1253.


 \bibitem [Sch] {Sch} Schaefer, H. H.
{\it Banach Lattices and Positive Operators}. Springer Verlag.
(1974).


 \bibitem [Sp] {Sp} Spitzer, F.
{\it A Variational characterization of finite Markov chains}. The
Annals of Mathematical Statistics. {\bf (43): N.1}  (1972), 303-307.

 \bibitem [Ta] {Ta} Takahashi, M. {\it Thermodynamics of one-dimensional solvable models}. Cambridge Press. 2005.

 \bibitem [Th] {Th} Thompson, C.
{\it Infinite-Spin Ising Model in one dimension}. Journal of
Mathematical Physics. {\bf (9): N.2} (1968), 241-245.


 \bibitem [V] {V}
van Enter, A., Romano, S. and  Zagrebnov, V. {\it First-order
transitions for some generalized $XY$ models}. J. Phys. A. {\bf 39},
no. {\bf 26}, (2006), 439-445.


 \bibitem [W] {W} Wrezinski, W. F. and Abdalla, E. {\it A precise formulation of the  third law of thermodynamics with applications to statistical physics and black holes}.
 Preprint USP (2007).


\end{thebibliography}
\end{document}